\documentclass[lecture]{pcms-l}
\makeatletter

\newif\iffirstlecture\firstlecturefalse

\newcommand{\lectureseries}{\firstlecturetrue
              \secdef\@lectureseries\@slectureseries} 

\newcommand{\@lectureseries}[2][default]{\chapter*{#2}%
              \gdef\thelectureseries{#1}} 

\newcommand{\@slectureseries}[1]{\chapter*{#1}}

\renewcommand{\auth}{\secdef\@auth\@sauth}

\newcommand{\@auth}[2][default]{\vspace{-1pc}{\raggedleft
        \Large\bf\noindent
        #2\endgraf
        \vspace*{2pc}
        }
        \def\@author{#1}%
}

\newcommand{\@sauth}[1]{\vspace{-1pc}{\raggedleft
        \Large\bf\noindent
        #1\endgraf
        \vspace*{2pc}
        }
        \def\@author{#1}%
}

\def\lecture#1{\global\Lecturetrue\global\Monographfalse
\iffirstlecture\else\chapter*{}\fi\firstlecturefalse
  \global\let\sectionmark\@gobble % \lecturemark will be used instead
  \addtocounter{lecture}1\relax
  \refstepcounter{chapter}%
\gdef\thelecturename{#1\unskip}
  {\Large\bfseries
   \raggedleft
   \@xp\uppercase\@xp{\thelecturelabel} {\LARGE\thelecturenum}\\
   \vspace*{3pt}%
   \thelecturename
   \endgraf}%
  \let\@secnumber=\thelecturenum
  \@xp\lecturemark\@xp{\thelecturename}%
  \addcontentsline{toc}{chapter}{%
    \thelecturelabel\ \thelecturenum.\ \thelecturename}%
  \vspace{34\p@}\noindent}
  
\def\lecture{\global\Lecturetrue\global\Monographfalse
\iffirstlecture\else\chapter*{}\fi%
  \global\let\sectionmark\@gobble % \lecturemark will be used instead
\secdef\@lecture\@slecture}

\def\@lecture[#1]#2{%
  \addtocounter{lecture}1\relax
  \refstepcounter{chapter}%
\gdef\thelecturename{#1\unskip}\firstlecturefalse
  {\Large\bfseries
   \raggedleft
   \@xp\uppercase\@xp{\thelecturelabel} {\LARGE\thelecturenum}\\
   \vspace*{3pt}%
    #2\unskip
   \endgraf}%
  \let\@secnumber=\thelecturenum
  \@xp\lecturemark\@xp{\thelecturename}%
  \addcontentsline{toc}{chapter}{%
    \thelecturelabel\ \thelecturenum.\ #2}%
  \vspace{34\p@}\noindent}
  
\def\slecturerunhead#1#2#3{%
    \let\@tempa\chaptername
    \uppercasenonmath{\@tempa}%
    \def\@tempb{#3\unskip}%
    \uppercasenonmath{\@tempb}%
    {\normalfont\@tempb}
    }
\def\slecturemark{%\let\@secnumber\@empty
    \@secmark\markright\slecturerunhead\chaptername}%

\def\@slecture#1{%
\iffirstlecture
\gdef\thelecturename{#1\unskip}\firstlecturefalse
  {\Large\bfseries
\noindent\thelecturename
   \endgraf}%
  \let\@secnumber=\thelecturenum
  \@xp\slecturemark\@xp{\thelecturename}%
  \addcontentsline{toc}{chapter}{%
    \thelecturename}%
 \vspace{-6\p@}\noindent
\else
\gdef\thelecturename{#1\unskip}\firstlecturefalse
  {\Large\bfseries
   \raggedleft
   \@xp\uppercase\@xp{\thelecturename}
   \endgraf}%
  \let\@secnumber=\thelecturenum
  \@xp\slecturemark\@xp{\thelecturename}%
  \addcontentsline{toc}{chapter}{%
    \thelecturename}%
  \vspace{34\p@}\noindent
\fi}

\ifLecture
  \def\chapterrunhead#1#2#3{%
    \let\@tempa\@author
    \uppercasenonmath{\@tempa}%
    \uppercasenonmath{\thelectureseries}%
    \textmd{\@tempa, \thelectureseries}
    }
  \def\lecturerunhead#1#2#3{%
    \let\@tempa\chaptername
    \uppercasenonmath{\@tempa}%
    \def\@tempb{#3\unskip}%
    \uppercasenonmath{\@tempb}%
    \textmd{\@tempa\ #2. \@tempb}
    }
\else
  \let\chapterrunhead\partrunhead
\fi

\newif\ifBibliographyIsASection\BibliographyIsASectionfalse

  \def\bibliomark{%\let\@secnumber\@empty
    \@secmark\markright\bibliorunhead\chaptername}%

  \def\bibliorunhead#1#2#3{%
    \let\@tempa\chaptername
    \uppercasenonmath{\@tempa}%
    \def\@tempb{#3\unskip}%
    \uppercasenonmath{\@tempb}%
    \textmd{\@tempb}
    }

\def\thebibliography#1{%
  \ifBibliographyIsASection
    \section*\refname
    \if@backmatter
      \markboth{\refname}{\refname}%
    \fi
  \else
\chapter*{}
  {\Large\bfseries
   \raggedleft
   \@xp\uppercase\@xp{\bibname} \\
   \endgraf}%
  \let\@secnumber=\thelecturenum
  \@xp\bibliomark\@xp{\bibname}%
  \addcontentsline{toc}{chapter}{%
    \bibname}%
  \vspace{34\p@}\noindent
  \fi
  \normalsize\labelsep .5em\relax
  \list{\@arabic\c@enumi.}{\settowidth\labelwidth{\@biblabel{#1}}%
  \leftmargin\labelwidth
  \advance\leftmargin\labelsep
	\usecounter{enumi}}\sloppy
  \clubpenalty9999 \widowpenalty\clubpenalty  \sfcode`\.\@m}

  \def\indexmark{%\let\@secnumber\@empty
    \@secmark\markright\indexrunhead\chaptername}%

  \def\indexrunhead#1#2#3{%
    \let\@tempa\chaptername
    \uppercasenonmath{\@tempa}%
    \def\@tempb{#3\unskip}%
    \uppercasenonmath{\@tempb}%
    \textmd{\@tempb}
    }

\ifLecture
\def\theindex{\cleardoublepage
\@restonecoltrue\if@twocolumn\@restonecolfalse\fi
\columnseprule \z@ \columnsep 35pt
\def\indexchap{\@startsection
		{chapter}{1}{\z@}{8pc}{34pt}%
		{\raggedleft
		\Large\bfseries}}%
 \twocolumn[\indexchap[{\indexname}]{\@xp\uppercase\@xp{\indexname}}]
  \@xp\indexmark\@xp{\indexname}%
	\thispagestyle{plain}\let\item\@idxitem\parindent\z@
	 \footnotesize\parskip\z@ plus .3pt\relax\let\item\@idxitem}
\fi

\def\@makefntext{\noindent\@makefnmark}

\def\setaddress{%
  {\let\@makefnmark\relax  \let\@thefnmark\relax
        \nobreak
        \addressnum@=\z@
        \loop\ifnum\addressnum@<\addresscount@\advance\addressnum@\@ne
           \footnote{$^{\hbox{\tiny\number\addressnum@}}$%
           \csname @address\number\addressnum@\endcsname
           \csname @curraddr\number\addressnum@\endcsname
           \csname @email\number\addressnum@\endcsname}\repeat
  \ifx\@empty\@date\else \@footnotetext{\@setdate}\fi
  \ifx\@empty\@subjclass\else \@footnotetext{\@setsubjclass}\fi
  \ifx\@empty\@keywords\else \@footnotetext{\@setkeywords}\fi
  \ifx\@empty\thankses\else \@footnotetext{%
    \def\par{\let\par\@par}\@setthanks}\fi
    }%
  \@setcopyright
}

\def\@tmpevenhead{\relax}

\def\cleardoublepage{\clearpage\if@twoside \ifodd\c@page\else
    \let\@tmpevenhead\@evenhead \let\@evenhead\relax\hbox{}\eject 
    \let\@evenhead\@tmpevenhead\if@twocolumn\hbox{}\newpage\fi\fi\fi}

\def\@setcopyright{%
  \let\copyrightyear\currentyear             % DF
  \insert\copyins{\hsize\textwidth
    \parfillskip\z@ \leftskip\z@\@plus.9\textwidth
    \fontsize{6}{7\p@}\normalfont\upshape
    \everypar{}%
    \vskip-\skip\copyins \nointerlineskip
    \noindent\vrule\@width\z@\@height\skip\copyins
    \copyright\copyrightyear\ American Mathematical Society\par
    \kern\z@}%
}

\renewcommand{\@auth}[2][default]{{\raggedleft
        \begingroup
  \fontsize{\@xivpt}{18}\bfseries%\centering
  #2\par \endgroup
        \vspace*{2pc}
        }
        \def\@author{#1}%
}

\renewcommand{\@sauth}[1]{{\raggedleft
        \begingroup
  \fontsize{\@xivpt}{18}\bfseries%\centering
  #1\par \endgroup
        \vspace*{2pc}
        }
        \def\@author{#1}%
}

\def\@lecture[#1]#2{%
  \addtocounter{lecture}1\relax
  \refstepcounter{chapter}%
\gdef\thelecturename{#1\unskip}\firstlecturefalse
  {\Large\bfseries
   \raggedleft
   \@xp\uppercase\@xp{\thelecturelabel} {\LARGE\thelecturenum}\\
   \vspace*{3pt}%
    #2\unskip
   \endgraf}%
  \let\@secnumber=\thelecturenum
  \@xp\lecturemark\@xp{\thelecturename}%
  \addcontentsline{toc}{chapter}{%
    \thelecturelabel\ \thelecturenum.\ #2}%
  \vspace{10\p@}\noindent}
  
\def\@slecture#1{%
\iffirstlecture
\gdef\thelecturename{#1\unskip}\firstlecturefalse
  {\Large\bfseries
\noindent\thelecturename
   \endgraf}%
  \let\@secnumber=\thelecturenum
  \@xp\slecturemark\@xp{\thelecturename}%
  \addcontentsline{toc}{chapter}{%
    \thelecturename}%
 \vspace{-6\p@}\noindent
\else
\gdef\thelecturename{#1\unskip}\firstlecturefalse
  {\Large\bfseries
   \raggedleft
   \@xp\uppercase\@xp{\thelecturename}
   \endgraf}%
  \let\@secnumber=\thelecturenum
  \@xp\slecturemark\@xp{\thelecturename}%
  \addcontentsline{toc}{chapter}{%
    \thelecturename}%
  \vspace{10\p@}\noindent
\fi}

\makeatother
\input epsf.tex

\usepackage{amssymb}

\def\currentvolume{??}
\def\currentyear{2003}

\numberwithin{section}{chapter}
\numberwithin{equation}{chapter}

\theoremstyle{plain}
\newtheorem{theorem}[equation]{Theorem}
\newtheorem{lemma}[equation]{Lemma}

\theoremstyle{definition}

\def\R{{\hbox{\bf R}}}
\def\C{{\hbox{\bf C}}}

\def\E{{\hbox{\bf E}}}
\def\P{{\hbox{\bf P}}}
\def\Z{{\hbox{\bf Z}}}
\def\T{{\hbox{\bf T}}}

\def\eps{{\varepsilon}}
\def\kp{{k'}}

\begin{document}

% if you do not want a title page and table of contents, omit
% the next three lines

\part*{Recent progress on the Restriction conjecture}
\pauth{Terence Tao}
\tableofcontents

\mainmatter
\setcounter{page}{1}

\LogoOn

\lectureseries[Recent progress on the Restriction conjecture]{Recent progress on the Restriction conjecture}

\auth{Terence Tao}

\address{Department of Mathematics, UCLA, Los Angeles CA 90095-1555}
\email{tao@math.ucla.edu}

%The following items will become first page footnotes; they are optional.

%\keywords{}
\date{June 30, 2003}
\thanks{The author is a Clay Prize Fellow and supported by a grant from the Packard foundation.  The author is also indebted to Jon Bueti, Tony Carbery, Andreas Seeger, and Jim Wright for several suggestions with the exposition and references, and especially to Julia Garibaldi for her careful reading of the text and exercises.}
\setaddress

\lecture{The restriction problem}

The purpose of these notes is describe the state of progress on the restriction problem in harmonic analysis, with an emphasis on the developments of the past decade or so on the Euclidean space version of these problems for spheres and other hypersurfaces.  As the field is quite large and has so many applications, it will be impossible to completely survey the field, but we will try to at least give the main ideas and developments in this area.

The restriction problem are connected to many other conjectures, notably the Kakeya and Bochner-Riesz conjectures, as well as PDE conjectures such as the local smoothing conjecture.  For reasons of space, we will not be able to discuss all these connections in detail; our main focus will be on proving restriction theorems for the Fourier transform.

Historically, the restriction problem originated by studying the Fourier transform of $L^p$ functions in Euclidean space $\R^n$ for some $n \geq 1$, although it was later realized that this problem also arises naturally in other contexts, such as non-linear PDE and in the study of eigenfunctions of the Laplacian.  

Fix $n \geq 2$; in our discussion all the constants $C$ are allowed to depend\footnote{The question on how quantifying the exact dependence of the constants here on the dimension $n$ as $n \to \infty$ is however an interesting problem, although to my knowledge there are not many results in this direction at present.} on $n$, and to vary from line to line.  

If $f$ is an $L^1(\R^n)$ function, then the Riemann-Lebesgue lemma implies that the Fourier transform $\hat f$, defined by $$ \hat f(\xi) := \int_{\R^n} e^{-2\pi i x \cdot \xi} f(x)\ dx$$ is a continuous bounded function on $\R^n$ which vanishes at infinity.  In particular, we can meaningfully restrict this function to any subset $S$ of $\R^n$, creating a continuous bounded function $\hat f|_S$ on $S$.

On the other hand, if $f$ is an arbitrary $L^2(\R^n)$ function, then the Fourier transform $\hat f$ can be any function in $L^2(\R^n)$, and in particular there is no meaningful way to restrict it to any set $S$ of zero measure.  

Between these two extremes, one may ask what happens to the Fourier transform of a function $f$ in $L^p(\R^n)$, where $1 < p < 2$.  Certainly we do not expect the Fourier transform $\hat f$ to be continuous or bounded, and it is easy to construct examples of $L^p$ functions which have an infinite Fourier transform at one point.  In fact, it is easy to create such a function which is infinite on an entire hyperplane; for instance, the function
\begin{equation}\label{ouch}
f(x) := \frac{\psi(x_2, \ldots, x_n)}{1 + |x_1|}
\end{equation}
where $x_1$ is the first co-ordinate of $x$ and $\psi$ is a bump function, lies in $L^p$ for every $p>1$, but has an infinite Fourier transform on every point on the hyperplane $\{ \xi \in \R^n: \xi_1 = 0 \}$.  One can of course concoct a similar example for any other hyperplane.

On the other hand, from the Hausdorff-Young inequality we see that $\hat f$ lies in the Lebesgue space $L^{p'}(\R^n)$, where $1/p + 1/p' = 1$.  Thus $\hat f$ can be meaningfully restricted to every set $S$ of positive measure.  

This leaves open the question of what happens to sets $S$ which have zero measure but which are not contained in hyperplanes.  In 1967 Stein made the surprising discovery that when such sets contain sufficient ``curvature'', that one can indeed restrict the Fourier transform of $L^p(\R^n)$ functions for certain $p > 1$.  This lead to the \emph{restriction problem} \cite{stein:problem}: for which sets $S \subseteq \R^n$ and which $1 \leq p \leq 2$ can the Fourier transform of an $L^p(\R^n)$ function be meaningfully restricted?  

There are of course infinitely many such sets to consider, but we shall focus our attention here on sets $S$ which are hypersurfaces\footnote{For surfaces of lower dimension, see \cite{christ:thesis}, \cite{prestini}, \cite{mock:habil}, \cite{banner}; for fractal sets in $\R$, see \cite{mock:fractal}, \cite{sj2}, \cite{mock:habil}; for surfaces in finite field geometries, see \cite{mock:tao}; for the restriction theory of the prime numbers, see \cite{green}.}, or compact subsets of hypersurfaces.  In particular, we shall be interested\footnote{It is easy to see, see Problem 1.1 below, using the symmetries of the Fourier transform, that the restriction problem for a set $S$ is unaffected by applying any translations or invertible linear transformations to the set $S$, so we can place the sphere, paraboloid, and cone in their standard forms \eqref{eq:sphere}, \eqref{eq:paraboloid}, \eqref{eq:cone} without loss of generality.} in the \emph{sphere}
\begin{equation}\label{eq:sphere}
S_{sphere} := \{ \xi \in \R^n: |\xi| = 1\},
\end{equation}
the \emph{paraboloid}
\begin{equation}\label{eq:paraboloid}
S_{parab} := \{ \xi \in \R^n: \xi_n = \frac{1}{2} |\underline \xi|^2 \},
\end{equation}
and the \emph{cone}
\begin{equation}\label{eq:cone}
S_{cone} := \{ \xi \in \R^n: \xi_n = |\underline \xi| \},
\end{equation}
where $\xi = (\underline{\xi}, \xi_n) \in \R^{n-1} \times \R \equiv \R^n$, and we always take $n \geq 2$ to avoid trivial situations.
These three surfaces are model examples of hypersurfaces with curvature\footnote{One could also consider cylinders such as $S^{k-1} \times \R^{n-k} \subset \R^n$, but it turns out that the restriction theory for these surfaces is identical to that of the sphere $S^{k-1}$ inside $\R^k$; see Problem 1.2 below.}, though of course the cone differs from the sphere and paraboloid in that it has one vanishing principal curvature.  These three hypersurfaces also enjoy a large group of symmetries (the orthogonal group, the parabolic scaling and Gallilean groups, and the Lorentz-Poincare group, respectively).  Also, these three hypersurfaces are related via the Fourier transform to solutions to certain familiar partial differential equations, namely the Helmholtz equation, Schr\"odinger equation, and wave equation; we will discuss this connection more in the last lecture.

\section{Restriction estimates: general theory}\label{sec:rest-extend}

Let $S$ be a compact subset (but with non-empty interior) of one of the above surfaces $S_{sphere}$, $S_{parab}$, $S_{cone}$.   We endow $S$ with a canonical measure $d\sigma$ - for the sphere, this is surface measure, for the parabola, it is the pullback of the $n-1$-dimensional Lebesgue measure $d\underline{\xi}$ under the projection map $\xi \mapsto \underline \xi$, while for the cone the pullback of $d\underline{\xi}/|\xi|$ is the most natural measure (as it is Lorentz-invariant; see Problem 1.2); thus
$$ \int_{S_{parab}} f(\xi) d\sigma(\xi) := \int_{\R^{n-1}} f(\underline{\xi}, \frac{1}{2} |\underline \xi|^2)\ d\underline \xi$$
and
$$ \int_{S_{cone}} f(\xi) d\sigma(\xi) := \int_{\R^{n-1}} f(\underline{\xi}, |\underline \xi|)\ \frac{d\underline \xi}{|\xi|}.$$
In order to restrict the Fourier transform of an $L^p(\R^n)$ function to $S$, it will suffice to prove an \emph{a priori} ``restriction estimate'' of the form\footnote{One may ask why we fixate on $L^p$ spaces (or more generally, Lorentz spaces such as $L^{p,\infty}$) here.  One reason is that these are the spaces which arise in Stein's maximal principle \cite{stein:maximal}; another is that these spaces are invariant under both translations and modulations.  Note that if $S$ is compact, then there is essentially no distinction between restricting the Fourier transform of an $L^p$ function $f$ and restricting a function in the Sobolev space $W^{s,p}$ for any $s \in \R$, since we only care about the frequencies of $f$ in a compact set.  It is however of interest to develop weighted estimates for restriction problems; there has been scant progress on this problem so far, but see \cite{stein:problem} for some conjectures, and \cite{carbery:weight} for some related work.}
\begin{equation}\label{eq:rpq}
 \| \hat f|_S \|_{L^q(S; d\sigma)} \leq C_{p,q,S} \| f \|_{L^p(\R^n)}
\end{equation}
for all Schwartz functions $f$ and some $1 \leq q \leq \infty$, since one can then use density arguments to obtain a continuous restriction operator from $L^p(\R^n)$ to $L^q(S; d\sigma)$ which extends the map $f \mapsto \hat f|_S$ for Schwartz functions.  When the set $S$ has sufficient symmetry (e.g. if $S$ is the sphere), this implication can in fact be reversed, using Stein's maximal principle \cite{stein:maximal}; if there is no bound of the form\footnote{Indeed, it suffices for the weak-type estimate from $L^p(\R^n)$ to $L^{p,\infty}(S; d\sigma)$ to fail.  See \cite{stein:maximal}; similar ideas arise in the factorization theory of Nikishin and Pisier.} \eqref{eq:rpq}, then one can construct functions $f \in L^p(\R^n)$ whose Fourier transform is infinite almost everywhere in $S$; see Problem 1.4 below.

The estimate \eqref{eq:rpq} can be written out more fully as
$$ (\int_S |\int_{\R^n} f(x) e^{-2\pi i x\cdot \xi}\ dx|^q d\sigma(\xi))^{1/q}
\leq C_{p,q,S} \|f\|_{L^p(\R^n)};$$
it is thus a model example of an \emph{oscillatory integral estimate}.  It is perhaps not surprising that this estimate is directly related to some other oscillatory integral estimates, and in particular the Bochner-Riesz and local smoothing estimates; more on this in a later lecture.

We will tend to think of $\R^n$ as representing ``physical space'', whose elements will be denoted names such as $x$ and $y$, while $S$ lives in ``frequency space'', and whose elements will be denoted names such as $\xi$ or $\omega$.  For the PDE applications it is sometimes convenient to think of $\R^n$ as a spacetime $\R^{n-1} \times \R := \{ (x,t): x \in \R^{n-1}, t \in \R \}$ (with the frequency space thus becoming spacetime frequency space $\{ (\xi,\tau): \xi \in \R^{n-1}, \tau \in \R \}$), but we will avoid doing so here.

It is thus of interest to see for which sets $S$ and which exponents $p$ and $q$ one has estimates of the form \eqref{eq:rpq}; henceforth we assume our functions $f$ to be Schwartz.  We denote\footnote{Strictly speaking, this should be $R_{S,d\sigma}(p \to q)$, since the choice of measure $d\sigma$ could be important; however in our contexts the measure $d\sigma$ will be clear from context, and in any event one can multiply $d\sigma$ by any measurable function bounded above and below without affecting the truth or falsity of $R_{S,d\sigma}(p \to q)$.} by $R_S(p \to q)$ the statement that \eqref{eq:rpq} holds for all $f$.    From our previous remarks we thus see that $R_S(1 \to q)$ holds for all $1 \leq q \leq \infty$, while $R_S(2 \to q)$ fails for all $1 \leq q \leq \infty$; the interesting question is then what happens for intermediate values of $p$.  If $S$ is compact, then an estimate of the form $R_S(p \to q)$ will automatically imply an estimate $R_S(\tilde p \to \tilde q)$ for all $\tilde p \leq p$ and $\tilde q \leq q$ by the Sobolev and H\"older inequalities (See Problem 1.5 below).  Thus the aim is to increase the size of $p$ and $q$ for which $R_S(p \to q)$ holds by as much as possible.

As mentioned before, if $d\sigma$ is Lebesgue measure then we have $R_S(p \to p')$ for all $1 \leq p \leq 2$ by Hausdorff-Young; if $S$ has finite measure then we indeed have $R_S(p \to q)$ for all $1 \leq p \leq 2$ and $q \leq p'$, by H\"older's inequality.  These are the only restriction estimates available in the finite measure case; see Problem 1.8.  The more interesting case is when $S$ has zero Lebesgue measure, and $d\sigma$ is a measure supported on $S$ (and thus singular to Lebesgue measure).

Also observe that if $S = S_1 \cup S_2$, then $R_S(p \to q)$ holds if and only if $R_{S_1}(p \to q)$ and $R_{S_2}(p \to q)$ both hold.  Thus the restriction property is a local property of the surface $S$; it does not depend, for instance, on the topology of $S$.  Also, from \eqref{ouch} we see that $R_S(p \to q)$ will fail for any $p > 1$ if $S$ contains a subset of a hyperplane which has positive measure (with respect to $\sigma$).

A simple duality argument shows that the estimate \eqref{eq:rpq} is equivalent to the ``extension estimate''
\begin{equation}\label{eq:rpq-dual}
 \| (F d\sigma)^\vee \|_{L^{p'}(\R^n)} \leq C_{p,q,S} \| F \|_{L^{q'}(S; d\sigma)}
\end{equation}
for all smooth functions $F$ on $S$, where $(F d\sigma)^\vee$ is the inverse Fourier transform of the measure $Fd\sigma$:
$$ (Fd\sigma)^\vee(x) := \int_S F(\xi) e^{2\pi i \xi \cdot x} d\sigma(\xi).$$
Indeed, the equivalence of \eqref{eq:rpq} and \eqref{eq:rpq-dual} follows from Parseval's identity 
$$\int_{\R^n} (F d\sigma)^\vee(x) \overline{f(x)}\ dx = \int_S F(\xi) \widehat{\overline{f(\xi)}}\ d\sigma(\xi)$$
and duality.  If we use $R^*_S(q' \to p')$ to denote the statement that the estimate \eqref{eq:rpq-dual} holds, then $R^*_S(q' \to p')$ is thus equivalent to $R_S(p \to q)$.

Note that because $F$ is smooth, it is possible to use the principle of stationary phase (see e.g. \cite{stein:large}) to obtain asymptotics for $(F d\sigma)^\vee$.  However, such asymptotics depend very much on the smooth norms of $F$, not just on the $L^{q'}(S)$ norm, and so do not imply estimates of the form \eqref{eq:rpq-dual} (although they can be used to provide counterexamples).  Thus one can think of extension estimates as a more general way than stationary phase to control oscillatory integrals, applicable in situations where the amplitude function $F(\xi)$ has magnitude bounds but no smoothness properties.

The extension formulation \eqref{eq:rpq-dual} also highlights the connection between this problem and partial differential equations.  For instance, consider a solution $u(t,x): \R \times \R^n \to \C$ to the free Schr\"odinger equation
$$ i\partial_t u + \Delta u = 0$$
with initial data $u(0,x) = u_0(x)$.  This has the explicit solution
$$ u(t,x) = \int e^{2\pi i (x \cdot \xi + 2\pi t |\xi|^2)} \hat u_0(\xi)\ d\xi,$$
or equivalently
$$ u = (F d\sigma)^\vee$$
where $d\sigma := d\xi \delta(\tau - 2\pi |\xi|^2)$ is (weighted) surface measure on the paraboloid $\{ (\tau,\xi) \in \R \times \R^n: \tau = 2 \pi |\xi|^2 \}$, and $F$ is the function $\hat u_0(\xi)$ restricted to the paraboloid.  Thus, estimates of the form $\R^*_S(q' \to p')$ when $S$ is the paraboloid in $\R \times \R^n$ to control certain spacetime norms of solutions to the free Schr\"odinger equation.  Somewhat similar connections exist between the cone \eqref{eq:cone} (in $\R \times \R^n$) and solutions to the wave equation $u_{tt} - \Delta u = 0$, or between the sphere \eqref{eq:sphere} and solutions to the Helmholtz equation $\Delta u + u = 0$.  We will return to these connections in a later lecture.

\section{Necessary conditions}\label{sec:necessary}

We will use the extension formulation \eqref{eq:rpq-dual} to develop some necessary conditions in order for $\R^*_S(q' \to p')$ to hold.  First of all, by setting $F \equiv 1$ we clearly see that we must have $(d\sigma)^\vee \in L^{p'}(\R^n)$ as a necessary condition.  In the case of the sphere \eqref{eq:sphere}, the Fourier transform $(d\sigma)^\vee(x)$ decays in magnitude like $(1 + |x|)^{-(n-1)/2}$ (as can be seen either by stationary phase, or by the asymptotics of Bessel functions), and so we obtain the necessary condition\footnote{There does not seem to be any hope for any weak-type endpoint estimate at $p' = 2n/(n-1)$, see \cite{bcss}.  One can also show that there are no estimates for $p' < 2n/(n-1)$ by using Problem 1.4 and the Knapp example.} $p' > 2n/(n-1)$, or equivalently $p < 2n/(n+1)$.  A similar computation gives the same constraint $p' > 2n/(n-1)$ for the paraboloid \eqref{eq:paraboloid}, while for the cone the asymptotics are slightly different, giving the condition $p' > 2(n-1)/(n-2)$.

Let $F$ be a smooth function on $S$ with an $L^\infty$ norm of at most 1.  Since $Fd\sigma$ is pointwise dominated by $d\sigma$, it seems intuitive that $(Fd\sigma)^\vee$ should be ``smaller'' than $(d\sigma)^\vee$.  Thus one should expect the above necessary conditions to in fact be sufficient to obtain the estimate $R^*_S(\infty \to p')$.  For completely general sets $S$, this assertion is essentially the \emph{Hardy-Littlewood majorant conjecture}; it is true when $p'$ is an even integer by direct calculation using Plancherel's theorem, but is false (for general sets $S$) for\footnote{A ``logarithmic'' failure was established by Bachelis in the 1970s; a more recent ``polynomial'' failure has been established independently by Mockenhaupt and Schlag (private communication) and Green and Ruzsa (private communication).  See \cite{mock:habil} for further discussion.} other values of $p'$.  However, it may still be that the majorant conjecture is still true for ``non-pathological'' sets $S$ such as the sphere, paraboloid, and cone.

Another necessary condition comes from the \emph{Knapp example} \cite{tomas:restrict}, \cite{strichartz:restrictionquadratic}.  In the case of the sphere or paraboloid, we sketch the example as follows.  Let $R \gg 1$.  Then, by a Taylor expansion of the surface $S$ around any interior point $\xi_0$, we see that the surface $S$ contains a ``cap'' $\kappa \subset S$ centered at $\xi_0$ of diameter\footnote{We use $X \lesssim Y$ or $X = O(Y)$ to denote an estimate of the form $X \leq CY$ where $C$ depends on $S$, $p$, $q$, but not on functions such as $f$, $F$, or on parameters such as $R$.  We use $X \sim Y$ to denote the estimate $X \lesssim Y \lesssim X$.} $\sim 1/R$ and surface measure $\sim R^{-(n-1)}$ which is contained inside a disk $D$ of radius $\sim 1/R$ and thickness $\sim 1/R^2$, oriented perpendicular to the unit normal of $S$ at $\xi_0$.  Let $F$ be the characteristic function of this cap $\kappa$ (one can smooth $F$ out if desired, but this does not affect the final necessary condition), and let $T$ be the dual tube to the disk $D$, i.e. a tube centered at the origin of length $\sim R^2$ and thickness $\sim R$ oriented in the direction of the unit normal to $S$ at $\xi_0$.  Then $(F d\sigma)^\vee$ has magnitude $\sim \sigma(K) \sim R^{-(n-1)}$ on a large portion of $T$ (this is basically because for a large portion of points $x$ in $T$, the phase function $e^{2\pi i x \cdot \xi}$ is essentially constant on $K$).  In particular, we have
$$ \| (F d\sigma)^\vee \|_{L^{p'}(\R^n)} \gtrsim |T|^{1/p'} R^{-(n-1)} \sim R^{-(n-1)} R^{(n+1)/p'},$$
while we have
$$ \| F \|_{L^{q'}(S; d\sigma)} \lesssim |\kappa|^{1/q'} \lesssim R^{-(n-1)/q'}.$$
Letting $R \to \infty$, we thus see that we need the necessary condition
$$ \frac{n+1}{p'} \leq \frac{n-1}{q}$$
in order for $R^*_S(q' \to p')$ to hold.  (In the case of compact subsets of the paraboloid with non-empty interior, one can obtain the same necessary condition using the parabolic scaling $(\underline \xi, \xi_n) \mapsto (\lambda \underline \xi, \lambda^2 \xi_n)$.  For the full (non-compact) paraboloid, one can improve this to $\frac{n+1}{p'} = \frac{n-1}{q}$; see Problem 1.6.).  In the case of the cone, we can lengthen the cap $\kappa$ in the null direction (so that it now has measure $\sim R^{-(n-2)}$ and lives in a ``plate'' of length $\sim 1$, width $\sim 1/R$ and thickness $1/R^2$), which eventually leads to the stronger necessary condition $\frac{n}{p'} \leq \frac{n-2}{q}$; as before, this can be strengthened to $\frac{n}{p'} = \frac{n-2}{q}$ if one is considering the full cone \eqref{eq:cone} and not just compact subsets of it with non-empty interior.

\begin{figure}[htbp] \centering
 \ \epsfbox{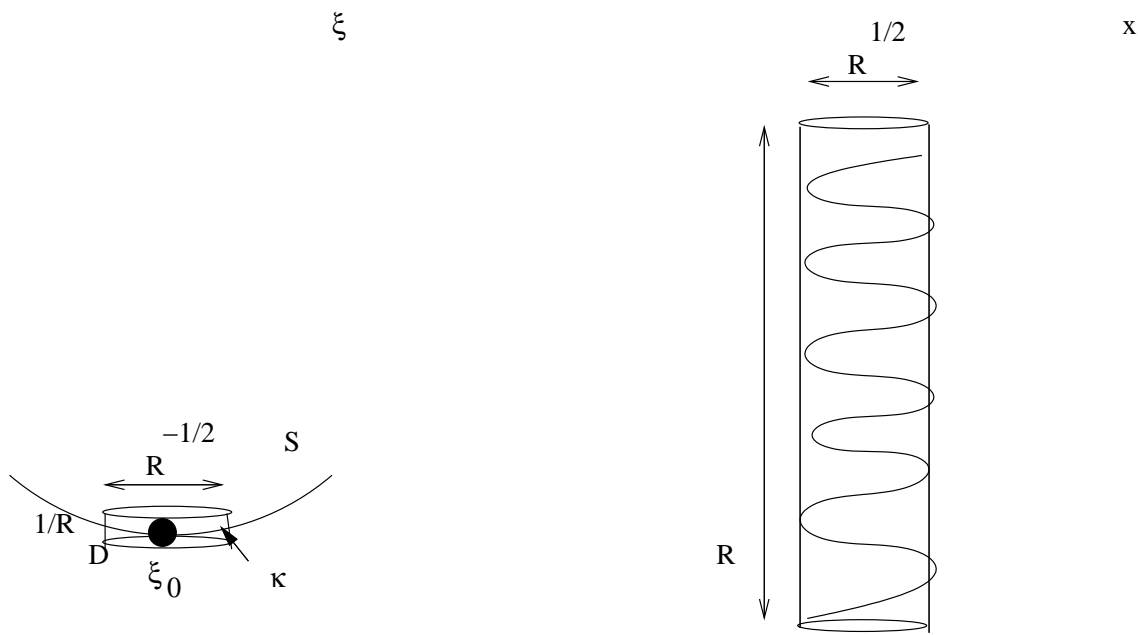}
 \caption{The (linear) Knapp example for a surface $S$ which has curvature at some point $\xi_0$.  The left
picture (labeled by $\xi$) represents frequency space; the function $F$ lives on the cap $\kappa$, which can
also be interpreted as the intersection of $S$ with a $R^{-1/2} \times R^{-1}$ disk.  The uncertainty principle then forces $(F d\sigma)^\vee$ in physical pace (labeled by $x$) to be mostly concentrated on a dual tube $T$ of dimensions $R^{1/2} \times R$ oriented in the direction normal to $D$.  The function $(F d\sigma)^\vee$ has some oscillation (basically of the form $e^{2\pi i x \cdot \xi}$ since $D$ is not centered at the origin; typically $|\xi| \sim 1$ so 
the oscillation has wavelength $1/|\xi| \sim 1$.  Note that $(Fd\sigma)^\vee$ will also be non-zero outside of $T$ but
there is usually some decay away from $T$ and so the portion inside $T$ is (heuristically at least) dominant.
}
\end{figure}

One can formulate a Knapp counterexample for any smooth hypersurface; the necessary conditions obtained this way become stronger as the surface becomes flatter, and in the extreme case where the surface is infinitely flat (e.g. when it is a hyperplane), there are no estimates.  See Problem 1.7.

\begin{figure}[htbp] \centering
 \ \epsfbox{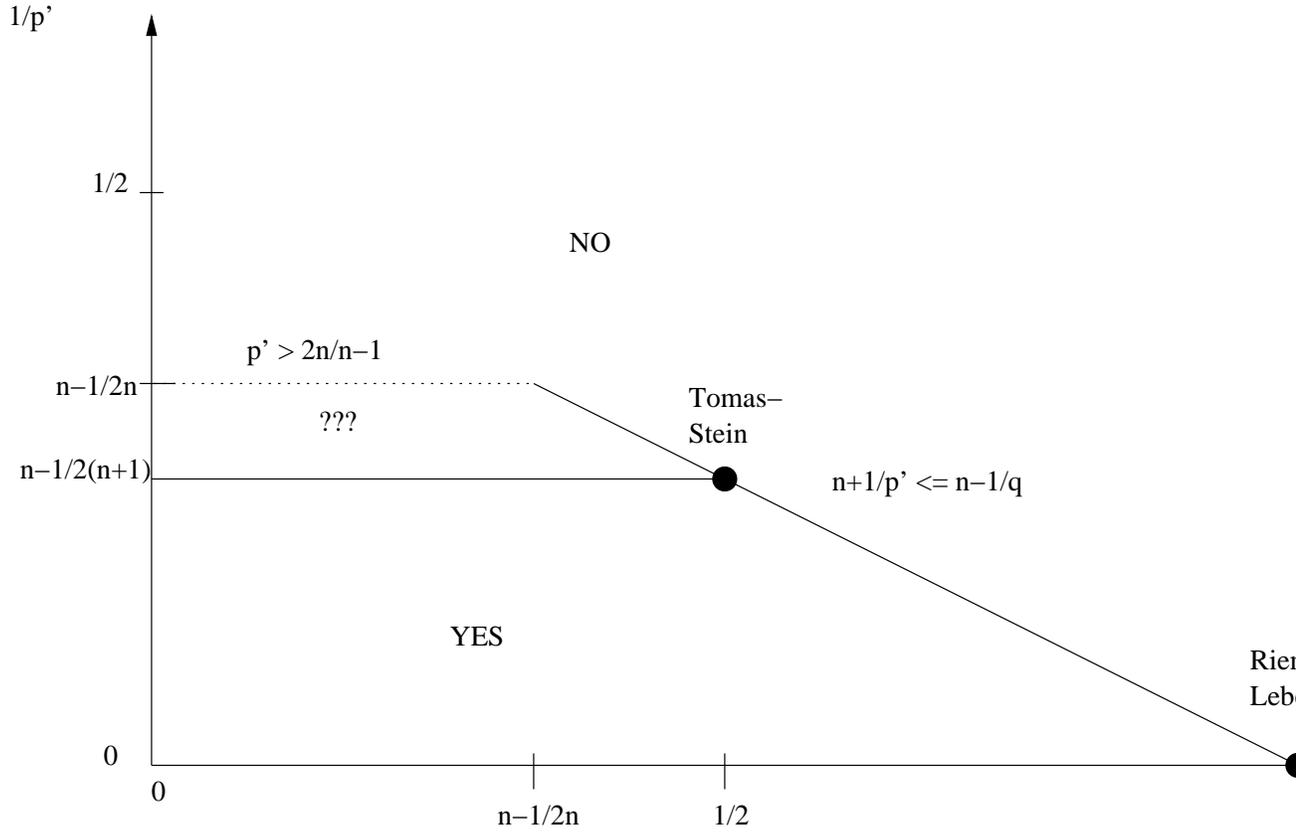}
 \caption{The range of exponents $p,q$ for the (linear) restriction problem for the sphere.  The necessary 
conditions $p' > 2n/(n-1)$ and $\frac{n+1}{p'} \leq \frac{n-1}{q}$ force $(1/q',1/p')$ to lie in the
trapezoidal region indicated.  The estimates on the bottom axis $1/p' = 0$ are very easy; the difficulty is
to make $1/p'$ as large as possible.  Once an estimate is obtained (e.g. the Tomas-Stein estimate $(1/q',1/p') = (1/2, \frac{n-1}{2(n+1)}$ displayed), one can use H\"older and interpolation to obtain all estimates to the left and below
of that estimate.  Thus to solve the restriction conjecture it will suffice to do so for $(1/q',1/p')$ arbitrarily
close to $(\frac{n-1}{2n}, \frac{n-1}{2n})$.
}
\end{figure}

The \emph{restriction conjecture} for the sphere, paraboloid, and cone then asserts that the above necessary conditions are in fact sufficient.  In other words, for compact subsets of the sphere and paraboloid the conjecture asserts that $R^*_S(q' \to p')$ holds when $p' > 2n/(n-1)$ and $\frac{n+1}{p'} \leq \frac{n-1}{q}$, while for compact subsets of the cone the conditions become $p' > 2(n-1)/(n-2)$ and $\frac{n}{p'} \leq \frac{n-2}{q}$ (i.e. the cone is conjectured to match the numerology of the sphere and paraboloid in one lower dimension; cf. Problem 1.3.).  This conjecture has been solved for the paraboloid and sphere in two dimensions, and for the cone in up to four dimensions; see Figures 1 and 2 for a more detailed summary of progress on this problem.  The restriction problems for the three surfaces are related; the sharp restriction conjecture for the sphere would imply the sharp restriction estimate for the paraboloid, because one can parabolically rescale the sphere to approach the paraboloid; see Problem 1.1 or \cite{tao:boch-rest}.  Also, using the method of descent, one can link the restriction conjecture for the cone in $\R^{n+1}$ with the restriction conjecture for the sphere, paraboloid, or other conic sections in $\R^n$, although the connection here is not as tight (see \cite{tao:cone} for some further discussion).

\begin{figure}\label{fig:parab}
\begin{tabular}{|l|l|l|} \hline
Dimension & Range of $p$ and $q$ & \\
\hline
$n = 2$ & $q' = 2, p' \geq 8$ & Stein, 1967 \\
& $q' > (p'/3)'; p' > 4$ & Fefferman, Stein, 1970 \cite{feff:thesis}\\
& $q' \geq (p'/3)'; p' > 4$ & Zygmund, 1974 \cite{zygmund} (best possible)\\
\hline
$n = 3$  & $q' = 2, p' \geq 6$ & Stein, 1967 \\

  & $q' > (p'/2)', p' > 4$ & Tomas 1975 \cite{tomas:restrict} \\
  & $q' \geq (p'/2)', p' \geq 4$ & Stein 1975; Sj\"olin $\sim$ 1975 \\
  & $q', p' > 4 - \frac{2}{15}$ & Bourgain 1991 \cite{borg:kakeya} \\
  & $q', p' > 4 - \frac{2}{11}$ & Wolff 1995 \cite{wolff:kakeya} \\
  & $q' > 7/3; p' > 4 - \frac{2}{11}$ & Moyua, Vargas, Vega 1996 \cite{vargas:restrict} \\
  & $q' \geq (p'/2)'; p' > 4 - \frac{5}{27}$ & Tao, Vargas, Vega 1998 \cite{tvv:bilinear}\\
  & $q' > 170/77; p' > 4 - \frac{2}{9}$ & Tao, Vargas, Vega 1998 \cite{tvv:bilinear}\\
  & $q' \geq (p'/2)'; p' > 4 - \frac{8}{31}$ & Tao, Vargas 2000 \cite{tv:cone1}\\
  & $q' > 26/11; p' > 4 - \frac{2}{7}$ & Tao, Vargas 2000 \cite{tv:cone1}\\
  & $q' \geq (p'/2)'; p' > 4 - \frac{2}{3}$ & Tao 2003 \cite{tao:parabola}\\
  & $q' \geq (p'/2)'; p' > 3$ & (conjectured)\\
\hline
$n > 3$ & $q' > ((n-1)p'/(n+1))'; p' > \frac{2(n+1)}{n-1}$ & Tomas 1975 \cite{tomas:restrict} \\
& $q' \geq ((n-1)p'/(n+1))'; p' \geq \frac{2(n+1)}{n-1}$ & Stein 1975 \\
& $q',p' > \frac{2(n+1)}{n-1} - \eps_n$ & Bourgain 1991 \cite{borg:kakeya} \\
& $q',p' > \frac{2n^2 + n + 6}{n^2 + n - 1}$  & Wolff 1995 \cite{wolff:kakeya} \\
& $q' > \frac{2n^2 + n + 6}{n^2 + 3n + 1}; p' > \frac{2n^2 + n + 6}{n^2 + n - 1}$  & Moyua, Vargas, Vega 1996 \cite{vargas:restrict} \\
  & $q' \geq ((n-1)p'/(n+1))'; p' > \frac{2(n+2)}{n}$ & Tao 2003 \cite{tao:parabola}\\
& $q' \geq ((n-1)p'/(n+1))'; p' > \frac{2n}{n-1}$ & (conjectured)\\
\hline
\end{tabular}
\caption{Known results on the restriction problem $R_S(p \to q)$ (or $R^*_S(q' \to p')$) for the sphere and for compact subsets of the paraboloid.  (For the whole paraboloid, restrict the above exponents to the range $q' = (\frac{(n-1)p'}{n+1})'$).}
\end{figure}

\begin{figure}\label{fig:cone}
\begin{tabular}{|l|l|l|} \hline
Dimension & Range of $p$ and $q$ & \\
\hline
$n = 3$  & $q' \geq (p'/3)', p' \geq 6$ & Strichartz 1977 \cite{strichartz:restrictionquadratic} \\
& $q' \geq (p'/3)'; p' > 4$ & Barcelo, 1985 \cite{barcelo} (best possible)\\
\hline
$n = 4$ &  $q' \geq (p'/2)', p' \geq 4$ & Strichartz 1977 
\cite{strichartz:restrictionquadratic} \\
& $q' \geq (p'/2)'; p' > 3$ & Wolff, 2000 \cite{wolff:cone} (best possible)\\
\hline
$n > 4$ &  $q' \geq ((n-2)p'/n)', p' \geq \frac{2n}{n-2}$ & Strichartz 1977 
\cite{strichartz:restrictionquadratic} \\
& $q' \geq ((n-2)p'/n)'; p' > \frac{2(n+2)}{n}$ & Wolff, 2000 \cite{wolff:cone}\\
& $q' \geq ((n-2)p'/n)'; p' > \frac{2(n-1)}{n-2}$ & (conjectured)\\
\hline
\end{tabular}
\caption{Known results on the restriction problem $R_S(p \to q)$ (or $R^*_S(q' \to p')$) for compact subsets of the cone.  (For the whole cone, restrict the above exponents to the range $q' = (\frac{(n-2)p'}{n})'$).}
\end{figure}

\section{Problems for Lecture 1}

\begin{itemize}

\item {\bf Problem 1.1.}  (a) Let $S$ be any subset of $\R^n$ with some measure $d\sigma$, and let $T: \R^n \to \R^n$ be any invertible affine transformation on $\R^n$ (i.e. $Tx := Lx + v$ for some fixed $v \in \R^n$ and some invertible linear transformation $L: \R^n \to \R^n$).  The image $T(S)$ of $S$ under the transform $T$ is thus endowed with the \emph{push-forward measure} $T_*(d\sigma)$, defined by
$$ \int_{T(S)} f(\xi)\ T_*(d\sigma)(\xi) := \int_S f(T \xi)\ d\sigma(\xi).$$
Show that for any $1\leq p,q \leq \infty$, the estimate $R_S(p \to q)$ holds if and only if $R_{T(S)}(p \to q)$ holds.
Furthermore, if $T$ is volume-preserving (i.e. $\det(L) = \pm 1$), then the best constant $C_{p,q,S}$ in \eqref{eq:rpq}
is the same for both $S$ and $T(S)$.  (Thus the restriction problem depends only on the shape of $S$, and not on its location or orientation).

\item (b)* Suppose that $1 \leq p, q \leq \infty$ obey the scaling relationship $\frac{n+1}{p'} = \frac{n-1}{q}$.  Suppose that the restriction estimate $R_{S_{sphere}}(p \to q)$ holds for the sphere.  Prove that the restriction
estimate $R_{S_{parab}}(p \to q)$ must then hold for the paraboloid.  (Hint: Translate the sphere $S_{sphere}$ upward by $e_n$ so that it touches the origin and is tangent to the hyperplane $\{ \xi_n = 0\}$.  Then apply a parabolic scaling $(\underline{\xi}, \xi_n) \mapsto (\lambda \underline{\xi}, \lambda^2 \xi_n)$ to the translated sphere, sending $\lambda \to \infty$.  As $\lambda \to \infty$, a Taylor expansion argument shows that this surface approaches the paraboloid $S_{parab}$.  The scaling condition will ensure that certain powers of $\lambda$ arising from Jacobians in the rescaling will cancel.  Now take limits, using for instance Fatou's lemma.)  As this result suggests, the restriction theory for the sphere and for the paraboloid are very closely related; indeed, in practice it has turned out that every result obtained for one has essentially also been obtained for the other.  See also \cite{carbery:parabola}, \cite{tao:boch-rest} for some similar connections.

\item {\bf Problem 1.2.}  Define the Minkowski form $\rho: \R^n \to \R$ by $\rho(\xi) := \xi_n^2 - |\underline \xi|^2$; thus $S_{cone}$ is the upper half of the zero set of $\rho$.  Show that the surface measure $d\sigma$ on the cone $S_{cone}$ is equal to the Dirac measure $\frac{1}{2} \delta(\rho(\xi)) H(\xi_n)$, where the Heaviside function $H(\xi_n)$ is equal to 1 when $\xi_n > 0$ and 0 when $\xi_n < 0$.  In other words, show that for any test function $F$ on $\R^n$, we have
$$ \int_{S_{cone}} F(\xi)\ d\sigma(\xi) = \lim_{\eps \to 0} \int_{\R^n} \frac{1}{2} \frac{\chi_{[-\eps,\eps]}(\rho(\xi))}{2\eps} H(\xi_n) F(\xi)\ d\xi.$$
Now let $L: \R^n \to \R^n$ be a Lorentz transformation (i.e. $L$ is linear and $\rho(L\xi) = \rho(\xi)$ for all $\xi \in\R^n$) which leaves the cone $S_{cone}$ invariant.  Show that $L$ also leaves $d\sigma$ invariant, i.e.
$$ \int_{S_{cone}} F(\xi)\ d\sigma(\xi) = \int_{S_{cone}} F(L(\xi))\ d\sigma(\xi)$$
for all test functions $F$.

\item This Lorentz invariance is very useful for the restriction theory of the cone; it allows one for instance to take a narrow sector of the cone and ``dilate'' it to a much wider sector, using a Lorentz boost in the direction of the sector.  See \cite{tv:cone2}, \cite{wolff:cone}, \cite{tao:cone}, \cite{wolff:smsub} for some applications of this technique.

\item {\bf Problem 1.3.}  Let $S_m$ be a compact subset of $\R^m$ with some finite non-zero measure $d\sigma_m$, and let $S_{n-m}$ be a compact subset of $\R^{n-m}$ with some finite non-zero measure $d\sigma_{n-m}$; the set $S_m \times S_{n-m} \subseteq \R^n$ is thus endowed with the product measure $d\sigma_m d\sigma_{n-m}$.  Let $1 \leq p,q \leq \infty$.  Show that the estimate $R_{S_m \times S_{n-m}}(p \to q)$ holds if and only if $R_{S_m}(p \to q)$ and $R_{S_{n-m}}(p \to q)$ both hold.  Conclude in particular that if $1 \leq p \leq 2$ and $p \leq q \leq p'$, that the restriction estimate $R_{S_{sphere}}(p \to q)$ holds for the sphere in $\R^n$ holds if and only if the restriction estimate $R_{S_{sphere} \times [0,1]}(p \to q)$ holds for the cylinder in $\R^{n+1}$.  Thus the restriction theory for the cylinder is essentially identical to that of the sphere of one lower dimension.  Note however that one cannot apply a similar argument to the cone.

\item {\bf Problem 1.4.} (a)  Let $\Omega$ be any measurable subset of the sphere $S_{sphere}$, which we endow with normalized surface measure $d\sigma$ (so $d\sigma(S_{sphere}) = 1$), and let $SO(n)$ be the group of rotations on $\R^n$, endowed with normalized Haar measure.  Let $R_1, R_2, \ldots, R_N$ be any $N$ rotations in $SO(n)$, chosen randomly and independently from $SO(n)$ using Haar measure as the probability measure.  The set $\bigcup_{j=1}^N R_j(\Omega)$ is then a subset of the sphere; show that the expected value of the $\sigma$-measure of this set is $1 - (1-|\Omega|)^N$.  Conclude in particular that if $|\Omega| \sim 1/N$, then there exists $N$ rotations $R_1, R_2, \ldots, R_N$ such that $\bigcup_{j=1}^N R_j(\Omega)$ has measure comparable to 1.

\item (b) Let $1 < p < \infty$ and $A > 0$.  Suppose that there exists a test function $f \in L^p(\R^n)$ and a $\lambda > 0$ such that
$$ d\sigma( \{ \xi \in S_{sphere}: |\hat f(\xi)| > \lambda \} ) \geq A^p \frac{\|f\|_{L^p(\R^n)}^p}{\lambda^p}.$$
Using (a), show that for any $0 < q < \infty$ there exists a function $F \in L^p(\R^n)$ such that
$$ \| \hat F|_{S_{sphere}} \|_{L^q(S_{sphere}; d\sigma)} \geq c_{p,q} A \| F \|_{L^p(\R^n)}$$ for some $c_{p,q} > 0$.  (Hint: We may assume that $f$ is compactly supported.  Set $\Omega := \{ \xi \in S_{sphere}: |\hat f(\xi)| > \lambda \}$, and apply (a) to cover a fair chunk of $S_{sphere}$ by rotations $R_j(\Omega)$ of $\Omega$.  Then set $F(x) := \sum_{j=1}^N \eps_j f(R_j(x - x_j))$, where the $\eps_j = \pm$ are randomized signs, and the $x_j$ are a sufficiently separated set of points.  Use Khinchin's inequality (see Appendix)).  

\item (c) Using (b), conclude that if the restriction operator $f \mapsto \hat f|_{S_{sphere}}$ is not of weak-type $(p,p)$, then it does not map $L^p(\R^n)$ to $L^q(S_{sphere}; d\sigma)$ for any $0 < q < \infty$.  (It is in fact possible to have a stronger conclusion - there exists a sequence $f^{(n)}$ of test functions converging in $L^p(\R^n)$ such that $\hat f^{(n)}|_S$ converges to infinity pointwise almost everywhere, but this is somewhat trickier to show).

\item {\bf Problem 1.5.}  Let $S$ be a compact subset of $\R^n$ with finite measure $d\sigma$, and let $1 \leq p_1 \leq p_2 \leq \infty$ and $1 \leq q_1 \leq q_2 \leq \infty$.  Show that if $R_S(p_2 \to q_2)$ holds, then $R_S(p_1 \to q_1)$ also holds.  (Hint: For the $q$ exponent, use H\"older.  For the $p$ exponent, use the fact that $\hat f|_S = \widehat{f*\phi}|_S$, where $\phi$ is any bump function whose Fourier transform equals 1 on a ball containing $S$).

\item {\bf Problem 1.6.}  Show that one can only obtain restriction estimates $R_K(p \to q)$ on any compact subset $K \subset S_{parab}$ of the paraboloid if $\frac{n+1}{p'} \leq \frac{n-1}{q}$, and one can only obtain restriction estimates $R_{S_{parab}}(p \to q)$ on the paraboloid if $\frac{n+1}{p'} = \frac{n-1}{q}$.  Also obtain a similar statement for the cone.

\item {\bf Problem 1.7.}  Let $k \geq 2$, and let $S$ be a surface of the form
$$ S := \{ (\underline{\xi}, \Phi(\underline{\xi})): \underline \xi \in \R^{n-1}, |\underline \xi| \leq 1 \}$$
where $\Phi: \R^{n-1} \to \R$ is a smooth function which vanishes to order $k$ at the origin, i.e. $\nabla^j \Phi(0) = 0$ for all $0 \leq j < k$.  We endow $S$ with the pull-back of Lebesgue measure $d\underline{\xi}$ under the projection map $(\underline{\xi}, \Phi(\underline{\xi})) \to \underline{\xi}$.  Show that one can only obtain restriction estimates $R_S(p \to q)$ when
$\frac{n+k-1}{p'} \leq \frac{n-1}{q}$.  In particular, if $\Phi$ is infinitely flat at the origin then we only have the trivial restriction estimates $R_S(1 \to q)$.

\item {\bf Problem 1.8.}  (a) Let $S = \R^n$ with Lebesgue measure.  Show that the restriction estimate $R_S(p \to q)$ holds if and only if $1 \leq p \leq 2$ and $q = p'$.  (Note that the ``if'' part is just the Hausdorff-Young inequality.  For the ``only if'' part, one can use some sort of scaling argument to obtain the $q=p'$ condition.  To obtain the $p \leq 2$ condition, try using a function $f$ which is a randomized sum of bump functions in different locations, using Khinchin's inequality).

\item (b)  Let $S$ be a subset of $\R^n$ with positive measure, endowed with Lebesgue measure.  Show that the restriction estimate $R_S(p \to q)$ holds if and only if $1 \leq p \leq 2$ and $q \leq p'$.

\end{itemize}

\section{Appendix: Khinchin's inequality}\label{sec:khinchin}

In this appendix we prove \emph{Khinchin's inequality}, which is fundamental in the use of randomization methods in harmonic analysis.  It concerns random sums of the form $\sum_{k=1}^N \eps_k f_k(x)$, where the $f_k(x)$ are functions and $\eps_k = \pm 1$ are random independent signs, with each $\eps_k$ equal to $+1$ with probability $1/2$ and $-1$ with probability $1/2$.  The intuition (coming from the law of large numbers) is that such a random sum should be distributed roughly like a Gaussian distribution around the origin with standard deviation $(\sum_{k=1}^N |f_k(x)|^2)^{1/2}$.  Khinchin's inequality is one way of making this intuition precise:

\begin{lemma}[Khinchin's inequality]  Let $0 < p < \infty$ and $f_1(x), \ldots, f_N(x)$ are a  collection of $L^p$-integrable functions on some measure space, and $\eps_1, \ldots, \eps_n$ are randomized
signs, then 
$$ \E(\| \sum_{k=1}^N \eps_k f_k \|_p^p) \sim 
\| (\sum_{k=1}^N |f_k|^2)^{1/2} \|_p^p,$$
where the constants in the $\sim$ symbol are independent of $N$ and
the $f_k$ (although they do depend on $p$), and $\E$ denotes
the expectation.
\end{lemma}

\begin{proof}
It suffices to show
\begin{equation}\label{khin}
 \E(| \sum_{k=1}^N \eps_k a_k |^p)^{1/p} \sim 
(\sum_{k=1}^N |a_k|^2)^{1/2},
\end{equation}
since one can then apply this inequality with $a_k = f_k(x)$ for each
$x$, raise this to the $p^{th}$ power, and integrate in $x$.  (Note
that the expectation operator $\E$ is linear).

It suffices to prove \eqref{khin} assuming that
$$ (\sum_{k=1}^N |a_k|^2)^{1/2} = 1.$$
This is because \eqref{khin} is unaffected by the operation of multiplying
$a_k$ by a constant.

We first prove this for $p=2$, in which case we have equality.  Indeed:

$$
\E(| \sum_{k=1}^N \eps_k a_k |^2) = \sum_{k=1}^N \sum_{\kp=1}^N
\E(\eps_k \eps_\kp a_k a_\kp)$$
by linearity of expectation.  By independence, the expectation
vanishes unless $k = \kp$, so our sum becomes
$$ \sum_{k=1}^N \E(\eps_k^2 a_k^2) = \sum_{k=1}^N a_k^2 = 1$$
as desired.

We have just proven that
$$  \E(| \sum_{k=1}^N \eps_k a_k |^2)^{1/2} = 1.$$

By H\"older's inequality this implies the upper bound
$$  \E(| \sum_{k=1}^N \eps_k a_k |^p)^{1/p} \leq 1$$
for all $0 < p \leq 2$.  We will now prove the upper bound
\begin{equation}\label{eq:khin-upper}
  \E(| \sum_{k=1}^N \eps_k a_k |^p)^{1/p} \leq C_p
(\sum_{k=1}^N |a_k|^2)^{1/2}
\end{equation}

for the remaining range $2 < p < \infty$.  The lower bound will then follow for all $0 < p < \infty$ since the quantity $\E(| \sum_{k=1}^N \eps_k a_k |^p)^{1/p}$ is log-convex in $p$ for all $0 < p < \infty$.

To show \eqref{eq:khin-upper} we first consider
the related expression
$$ \E(e^{\lambda \sum_{k=1}^N \eps_k a_k})$$
where $\lambda > 0$ is a parameter.  This is of course equal to
$$ \E(\prod_{k=1}^N e^{\lambda \eps_k a_k}).$$
By independence one can take the product outside of the expectation:
$$ \prod_{k=1}^N \E(e^{\lambda \eps_k a_k}) = \prod_{k=1}^N \cosh(\lambda a_k).$$
By comparing Taylor series, we see that $\cosh(x) \leq e^{x^2/2}$ for all
$x$, so we have
$$ \E(e^{\lambda \sum_{k=1}^N \eps_k a_k}) \leq \prod_{k=1}^N e^{\lambda^2
a_k^2/2}
= e^{\lambda^2/2}.$$
From the Chebyshev inequality
$$ \P( X > \alpha ) \leq \frac{1}{\alpha} \E(X)$$
where $\P(E)$ is the probability of an event $E$, we obtain
$$ \P( e^{\lambda \sum_{k=1}^N \eps_k a_k} > \alpha )
\lesssim \alpha^{-1} e^{\lambda^2/2}.$$
This works for every $\alpha$; we choose $\alpha = e^{\lambda^2}$.
The estimate now becomes
$$ \P( \sum_{k=1}^N \eps_k a_k > \lambda ) \lesssim e^{-\lambda^2/2}.$$
Since the random variable $\sum_{k=1}^N \eps_k a_k$ is clearly
symmetric around the origin, we thus have
$$ \P( |\sum_{k=1}^N \eps_k a_k| > \lambda ) \lesssim e^{-\lambda^2/2}.$$
In particular, we have
$$ \P( |\sum_{k=1}^N \eps_k a_k|^p > \lambda ) \lesssim e^{-\lambda^{2/p}/2}.$$
If we now use the identity
$$ \E(X) = \int_0^\infty \P(X > \lambda)\ d\lambda$$
for any random variable $X$ with finite expectation (which follows immediately from Fubini's theorem), we thus
obtain
$$ \E(X) \lesssim 1$$
for all $1 < p < \infty$, as desired.
\end{proof}

\lecture{Some tools used to prove restriction estimates.}

We now begin discussing some of the tools used to prove the above restriction theorems.  In this lecture we discuss two basic ones: firstly, the decay of the Fourier transform of surface measure and how it can be used to localize restriction estimates; and secondly, the bilinear approach to restriction estimates which are especially good for $L^4$ estimates, but are also very useful in the $L^p$ theory as they mostly eliminate the problem of small angles.

There are two more techniques, both very powerful, which are also needed to prove the best restriction estimates\footnote{In fact, the most recent restriction estimate - a sharp bilinear restriction estimate for paraboloids \cite{tao:parabola} - uses virtually every single technique mentioned in these lectures!} known to date: the \emph{wave packet decomposition}, and the \emph{induction on scales} method.  These however are more complicated and will be discussed in later lectures.

\section{Local restriction estimates}\label{sec:local-rest}

The first key idea is to reduce the study of \emph{global} restriction theorems (where the physical space variable is allowed to range over all of $\R^n$), to that of \emph{local} restriction theorems (where the physical space variable is constrained to lie in a ball).  As we shall see, the reason we can obtain this reduction is because the Fourier transform $(d\sigma)^\vee(x)$ of surface measure $d\sigma$ is somewhat localized in space (i.e. it decays as $|x| \to +\infty$).

More precisely, for any exponents $p,q$, and any $\alpha \geq 0$, let $R_S(p \to q; \alpha)$ denote the statement that the localized restriction estimate
\begin{equation}\label{eq:local-rest}
  \| \hat f|_S \|_{L^q(S; d\sigma)} \leq C_{p,q,S,\alpha} R^\alpha \| f \|_{L^p(B(x_0,R))}
\end{equation}
holds for any radius $R \geq 1$, any ball $B(x_0,R) := \{ x \in \R^n: |x-x_0| \leq R\}$ of radius $R$, and any test function $f$ supported in $B(x_0,R)$.  Note that the center $x_0$ of the ball is irrelevant since one can translate $f$ by an arbitrary amount without affecting the \emph{magnitude} of $\hat f$.  The condition $\alpha \geq 0$ is necessary since the claim \eqref{eq:local-rest} is clearly absurd for negative $\alpha$ (at least if $\sigma$ has non-zero total measure), as can be seen by letting $R \to \infty$.

Observe that estimates for lower $\alpha$ immediately imply estimates for higher $\alpha$ (keeping $p$, $q$, $S$ fixed).  Also, the local estimate $R_S(p \to q; 0)$ is clearly equivalent to the global estimate $R_S(p \to q)$ by a sending $R \to \infty$ and applying a limiting argument.  Finally, it is easy to prove estimates of this type for very large $\alpha$; for instance, for smooth compact hypersurfaces $S$ one has the estimate $R_S(p \to q; n/p')$ just from the H\"older inequality
$$ \|\hat f \|_{L^{\infty}(S; d\sigma)} \leq \| f\|_{L^1(B(x_0, R))} \leq C_p R^{n/p'} \| f \|_{L^p(B(x_0,R))}.$$
Thus the aim is to lower the value of $\alpha$ from the trivial value of $\alpha = n/p'$, toward the ultimate aim of $\alpha = 0$, at least when $p$ and $q$ lie inside the conjectured range of the restriction conjecture.  (For other $p$ and $q$, the canonical counterexamples will give some non-zero lower bound on $\alpha$; see Problem 2.1.).

By duality, the local restriction estimate $R_S(p \to q; \alpha)$ is equivalent to the local extension estimate $R^*_S(q' \to p'; \alpha)$, which asserts that
\begin{equation}\label{eq:local-extension}
\| (F d\sigma)^\vee \|_{L^{p'}(B(x_0,R))} \leq C_{p,q,S,\alpha} R^\alpha \| F \|_{L^{q'}(S; d\sigma)}
\end{equation}
for all smooth functions $F$ on $S$, all $R \geq 1$, and all balls $B(x_0,R)$.

The uncertainty principle suggests that since the spatial variable has now been localized to scale $R$, the frequency variable can be safely blurred to scale $1/R$.  In the case where $S$ is a smooth compact hypersurface, this is indeed correct; the estimate \eqref{eq:local-rest} is equivalent to the estimate
\begin{equation}\label{eq:local-rest-blur}
  \| \hat f \|_{L^q(N_{1/R}(S))} \leq C_{p,q,S,\alpha} R^{\alpha-1/q} \| f \|_{L^p(B(x_0,R))}
\end{equation}
holding for all test functions $f$ supported on $B(x_0,R)$, where $N_{1/R}(S)$ is the $1/R$-neighborhood of $S$; see Problem 2.2 for this fact and a slight refinement.

Of course, \eqref{eq:local-rest-blur} is equivalent by duality to the estimate
\begin{equation}\label{eq:blur-adjoint}
\| G^\vee \|_{L^{p'}(B(x_0,R))} \leq C_{p,q,S,\alpha} R^{\alpha-1/q} \| G \|_{L^{q'}(N_{1/R}(S))}
\end{equation}
for all smooth functions $G$ supported on $N_{1/R}(S)$.  By another application of the uncertainty principle (similar to the transference principle of Marcinkiewicz and Zygmund; see Problem 2.3), this estimate is also equivalent (when $q \leq p'$) to the discrete version\footnote{Thus restriction estimates give some information on the magnitude of exponential sums when one only has size information on the coefficients $g(\xi)$ and not phase information.  One can also connect restriction theorems to the related topic of $\Lambda(p)$ sets; see e.g. \cite{mock:fractal} for a discussion.}
\begin{equation}\label{eq:discrete-adjoint}
\| \sum_{\xi \in \Lambda} g(\xi) e^{2\pi i x \cdot \xi} \|_{L^{p'}(B(x_0,R))} \leq C_{p,q,S,\alpha} R^{\alpha+(n-1)/q} \| g \|_{l^{q'}(\Lambda)}
\end{equation}
where $\Lambda$ is any maximal $1/R$-separated subset of $S$, and $g: \Lambda \to \C$ is any (discrete) function on $\Lambda$.

From the formulation \eqref{eq:local-rest-blur} and Plancherel's theorem, we immediately obtain the local restriction estimate $R_S(2 \to 2; 1/2)$ for smooth compact hypersurfaces $S$; this estimate can also be obtained from the Agmon-H\"ormander estimate or from the frequency-localized version of the Sobolev trace lemma.

To convert local restriction estimates into global ones, the key tool used is the decay of the Fourier transform $(d\sigma)^\vee$.  Indeed, suppose we have a decay estimate of the form
\begin{equation}\label{eq:decay-ass}
 |(d\sigma)^{\vee}(x)| \leq \frac{C}{(1+|x|)^\rho}
\end{equation}
for some $\rho > 0$.  (For the sphere and compact subsets of the paraboloid, this estimate holds for $\rho = (n-1)/2$; for the cone, this holds for $\rho = (n-2)/2$.  See Problem 2.4.)  Then the contributions to \eqref{eq:rpq} arising from widely separated portions of space will be almost orthogonal.  For instance, suppose $R \geq 1$ and $B(x_0,R)$ and $B(x_1,R)$ are two balls which are separated by at least a distance of $R$.  Then if $f_0$ and $f_1$ are supported on $B(x_0,R)$ and $B(x_1,R)$ respectively, the Fourier transforms $\hat f_0|_S$ and $\hat f_1|_S$ will be almost orthogonal on $S$:
\begin{equation}\label{eq:decay}
\begin{split}
 |\langle \hat f_0|_S, \hat f_1|_S \rangle_{L^2(S; d\sigma)}| &= |\langle \hat f_0 d\sigma, \hat f_1 \rangle_{L^2(\R^n)}|\\
 &= |\langle f_0 * (d\sigma)^\vee, f_1 \rangle|\\
 &\leq C R^{-\rho} \|f_0\|_{L^1(B(x_0,R))} \| f_1\|_{L^1(B(x_1,R))},
\end{split}
\end{equation}
since the convolution kernel $(d\sigma)^\vee$ has magnitude $O(R^{-\rho})$ when applied to differences of points in $B(x_0,R)$ and points in $B(x_1,R)$.  This almost orthogonality asserts in some sense that distant balls do not interact much with each other, and so will allow us to reduce a global restriction estimate to a local one.

One heuristic way to view \eqref{eq:decay} is as follows.  This estimate is in some sense a ``bilinear'' version of the (false) estimate
\begin{equation}\label{eq:false}
 \| \hat f_0 \|_{L^2(S; d\sigma)} \leq C R^{-\rho/2} \| f_0 \|_{L^1(B(x_0,R))};
\end{equation}
this estimate is of course not true since the limit of the estimate as $R \to \infty$ is absurd, nevertheless it is ``virtually'' true in the sense that it implies the true estimate \eqref{eq:decay} by Cauchy-Schwarz.  Note that \eqref{eq:false} is just the (false) local restriction estimate $R_S(1 \to 2; -\rho/2)$.  While this estimate is not true, it is true for certain interpolation purposes; for instance, by combining it with the Agmon-H\"ormander estimate $R_S(2 \to 2; 1/2)$, one can obtain the \emph{Tomas-Stein estimate} $R_S(\frac{2(\rho+1)}{\rho+2} \to 2)$, or more generally 
\begin{equation}\label{ts}
R_S(p \to 2) \hbox{ for all } p \leq \frac{2(\rho+1)}{\rho+2}.
\end{equation}
This heuristic argument can be made rigorous by using orthogonality arguments such as the $TT^*$ method; see Problem 2.5.  In the particular cases of the sphere and paraboloid, the Tomas-Stein estimate yields $R_S(\frac{2(n+1)}{n+3} \to 2)$; for the cone, it yields $R_S(\frac{2n}{n-2} \to 2)$.  Note that this is consistent with the numerology supplied by the Knapp example from the previous lecture.

The Tomas-Stein argument uses orthogonality on $L^2(S; d\sigma)$, and at first glance it thus seems that it can only be applied to obtain restriction theorems $R_S(p \to q)$ when $q=2$.  However, it was observed by Bourgain \cite{borg:kakeya}, \cite{borg:stein} that the same type of orthogonality arguments, exploiting the decay of the Fourier transform of $d\sigma$, can also be used for restriction theorems which are not $L^2$-based, albeit with some inefficiencies due to the use of non-$L^2$ orthogonality estimates.  These ideas were then extended in \cite{vargas:restrict}, \cite{borg:cone}, \cite{tao:boch-rest}, \cite{tvv:bilinear}, \cite{tv:cone1}; we cite two sample results (of a rather technical nature) below.

\begin{theorem}\label{borg:global}\cite{borg:kakeya}, \cite{borg:stein}, \cite{vargas:restrict}, \cite{tvv:bilinear}, \cite{tv:cone1} Let $\rho$ be as above.  If $R^*_S(p \to q; \alpha)$ holds for some $\rho + 1 > \alpha q$, then we have $R^*(\tilde p \to \tilde q)$ whenever
$$ \tilde q > 2 + \frac{q}{\rho + 1 - \alpha q}; \quad \frac{\tilde p}{\tilde q} < 1 + \frac{q}{p(\rho + 1 - \alpha q)}.
$$
\end{theorem}

\begin{theorem}\label{thm:eps-remove}\cite{tao:boch-rest}, \cite{tao:weak2} Let $\rho$ be above.  If $R_S(p \to p; \alpha)$ holds for some $p < 2$ and $0 < \alpha \ll 1$, then we have $R_S(p \to q)$ whenever
$$ \frac{1}{q} > \frac{1}{p} + \frac{C_\rho}{\log(1/\alpha)}.$$
\end{theorem}

The second theorem in particular has the following consequence: if $R_S(p \to p; \eps)$ is true for all $\eps > 0$, then $R_S(p \to p-\eps)$ is also true for every $\eps > 0$.  (The converse statement follows easily from interpolation).  Thus we can convert a local estimate with epsilon losses to a global estimate, where the epsilon loss has now been transferred to the exponents.  This type of ``epsilon-removal lemma'' is common in this theory, see \cite{borg:cone}, \cite{tv:cone1}, \cite{tao:weak2} (or Theorem \ref{thm:eps} below) for some more examples.

The above results are probably not optimal, however they do emphasize the point that one can study global restriction estimates via their local counterparts.

\section{Bilinear restriction estimates}\label{sec:bilinear}

We now turn to another idea in the development of restriction theory - that of passing from the linear restriction and extension estimates to bilinear analogues.

The original motivation of this theory was the ``$L^4$'' or ``bi-orthogonality'' theory developed in such places as \cite{feff:note}, \cite{cordoba:sieve}, \cite{carl:disc}, \cite{carbery:maximal-bochner}, \cite{mock:cone}.  The basic idea is that expressions such as $\| (F d\sigma)^\vee \|_{L^{p'}(\R^n)}$ can be calculated very explicitly when $p'$ is an even integer, and in particular when $p'$ is equal to 4.  Indeed, we have by Plancherel's theorem that
$$ \| (F d\sigma)^\vee \|_{L^4(\R^n)} = \| (F d\sigma)^\vee (F d\sigma)^\vee \|_{L^2(\R^n)}^{1/2} = \| F d\sigma * F d\sigma \|_{L^2(\R^n)}^{1/2}.$$
Thus one can reduce a restriction estimate such as $R_S^*(q' \to 4)$ to an estimate of the form
$$ \| F d\sigma * F d\sigma \|_{L^2(\R^n)} \leq C_q \| F \|_{L^{q'}(S; d\sigma)}^2;$$
the point here is that there is no oscillation in this estimate (since there is no Fourier transform), and this estimate can be proven or disproven by more direct methods.  For instance when $S$ is the circle in $\R^2$, there is a logarithmic divergence in the above estimate, since $d\sigma * d\sigma$ blows up like $1/|x|^{1/2}$ on the circle $\{ x \in \R^2: |x| = 2\}$ of radius 2, however by introducing the localizing parameter $R$ one can easily prove the modified estimate
\begin{equation}\label{eq:gg}
 \| G * G \|_{L^2(\R^n)} \lesssim (\log R)^{1/2} R^{-3/2} \| G \|_{L^4(N_{1/R}(S))}^2,
\end{equation}
for all $R \geq 1$ and all $G$ supported on $N_{1/R}(S)$ (Problem 2.6); comparing this with \eqref{eq:blur-adjoint} we obtain the local restriction estimate $R_S^*(4 \to 4, \eps)$ for any $\eps > 0$, which (by use of epsilon-removal lemmas such as Theorem \ref{thm:eps-remove}) proves the optimal range of restriction estimates for the circle (first proven by Zygmund \cite{zygmund}, by a more direct argument).

Similar arguments also give the optimal restriction theory for the cone in three dimensions, see \cite{barcelo}.  At first glance, this theory seems to be limited to $L^4$ (or to $L^6$, $L^8$ etc.), since it relies on Plancherel's theorem.  However, one can partially extend these ideas to other exponents $L^{p'}$, even when $p'$ is not an even integer.  The main point is that the linear estimate
$$ \| (F d\sigma)^\vee \|_{L^{p'}(\R^n)} \leq C_{p,q,S} \| F \|_{L^{q'}(S; d\sigma)}$$
is equivalent, via squaring, to the quadratic estimate
$$ \| (F d\sigma)^\vee (F d\sigma)^\vee \|_{L^{p'/2}(\R^n)} \leq C_{p,q,S} \| F \|_{L^{q'}(S; d\sigma)} \| F \|_{L^{q'}(S; d\sigma)}$$
which we can depolarize as the bilinear estimate
\begin{equation}\label{eq:bilinear}
 \| (F_1 d\sigma)^\vee (F_2 d\sigma)^\vee \|_{L^{p'/2}(\R^n)} \leq C_{p,q,S} \| F_1 \|_{L^{q'}(S; d\sigma)} \| F_2 \|_{L^{q'}(S; d\sigma)}.
\end{equation}
In such an estimate, the worst case typically occurs when $F_1$ and $F_2$ are both concentrated in the same small ``cap'' in $S$; this is what happens in the Knapp example, for instance.  

The strategy of the bilinear approach to restriction theory is to rewrite the linear estimate \eqref{eq:rpq-dual} as the bilinear estimate \eqref{eq:bilinear}, which in turn is a special case of a more general estimate of the form
\begin{equation}\label{eq:bilinear-gen}
\| (F_1 d\sigma_1)^\vee (F_2 d\sigma_2)^\vee \|_{L^{p'/2}(\R^n)} \leq C_{p,q, S_1, S_2} \| F_1 \|_{L^{q'}(S_1; d\sigma_1)} \| F_2 \|_{L^{q'}(S_2; d\sigma_2)},
\end{equation}
for arbitrary \emph{pairs} of smooth compact hypersurfaces $S_1$, $S_2$ with surface measures $d\sigma_1$, $d\sigma_2$ respectively, and all smooth $F_1$, $F_2$ supported on $S_1$ and $S_2$.  We let $R^*_{S_1,S_2}(q' \times q' \to p'/2)$ denote the statement that the estimate \eqref{eq:bilinear-gen} holds.  Then by the above discussion, $R^*_S(q' \to p')$ is equivalent to $R^*_{S,S}(q' \times q' \to p'/2)$.  Thus linear restriction estimates are special cases of bilinear estimates.  However, there are bilinear estimates that cannot be derived directly from linear ones.  For instance, let $S_1 := \{ (\xi_1,0): \xi_1 \in \R \}$ and $S_2 := \{ (0,\xi_2): \xi_2 \in \R \}$ denote the $x$ and $y$ axes in $\R^2$.  Then we have $(F_1 d\sigma_1)^\vee(x,y) = \check F_1(x)$ and $(F_2 d\sigma_2)^\vee(x,y) = \check F_2(y)$, and so there are no global restriction estimates of the form $R^*_{S_1}(q' \to p')$ or $R^*_{S_2}(q' \to p')$ unless $p' = \infty$, since the Fourier transforms do not decay at infinity.  However, since
$$(F_1 d\sigma_1)^\vee (F_2 d\sigma_2)^\vee(x,y) = \check F_1(x) \check F_2(y),$$
we see from the one-dimensional Plancherel theorem that we have the bilinear restriction estimate $R^*_{S_1, S_2}(2 \times 2 \to 2)$.  Note however that the symmetrized analogues $R^*_{S_1,S_1}(2 \times 2 \to 2)$ and $R^*_{S_2,S_2}(2 \times 2 \to 2)$ are false.  Thus the bilinear estimate exploits the \emph{transversality} of $S_1$ and $S_2$.

A higher-dimensional analogue of this estimate is known: if $S_1$ and $S_2$ are two smooth compact hypersurfaces which are \emph{transverse} in the sense that the set of unit normals of $S_1$ are separated by some non-zero distance from the set of unit normals of $S_2$, then we have $R^*_{S_1, S_2}(2 \times 2 \to 2)$.  This can be easily seen by using Plancherel to convert the bilinear restriction estimate to a bilinear convolution estimate
$$
\| (F_1 d\sigma_1) * (F_2 d\sigma_2) \|_{L^{2}(\R^n)} \leq C_{p,q, S_1, S_2} \| F_1 \|_{L^{2}(S_1; d\sigma_1)} \| F_2 \|_{L^{2}(S_2; d\sigma_2)},$$
and then using the Cauchy-Schwarz estimate
$$
\| (F_1 d\sigma_1) * (F_2 d\sigma_2) \|_{L^{2}(\R^n)}
\leq \| (|F_1|^2 d\sigma_1) * (|F_2|^2 d\sigma_2) \|_{L^1(\R^n)}
\| d\sigma_1 * d\sigma_2 \|_{L^\infty(\R^n)}$$
and using transversality to bound the second factor (this is a generalization of the argument in Problem 2.6).  Generalizations of these ``bilinear $L^2$'' estimates have arisen in recent work in non-linear evolution equations (starting with the work of Bourgain \cite{borg:xsb} and Klainerman-Machedon \cite{klainerman:nulllocal} and continued by numerous authors, see for instance \cite{kpv:kdv}) and are especially useful for handling non-linearities which contain derivatives arranged to create a ``null form'', but we will not pursue this matter here, and refer the reader instead to \cite{ginibre:survey}, \cite{damiano:null}, \cite{tao:xsb}.  There has been also some work in generalizing these bilinear estimates to weighted settings, see \cite{bcc}.

\begin{figure}[htbp] \centering
 \ \epsfbox{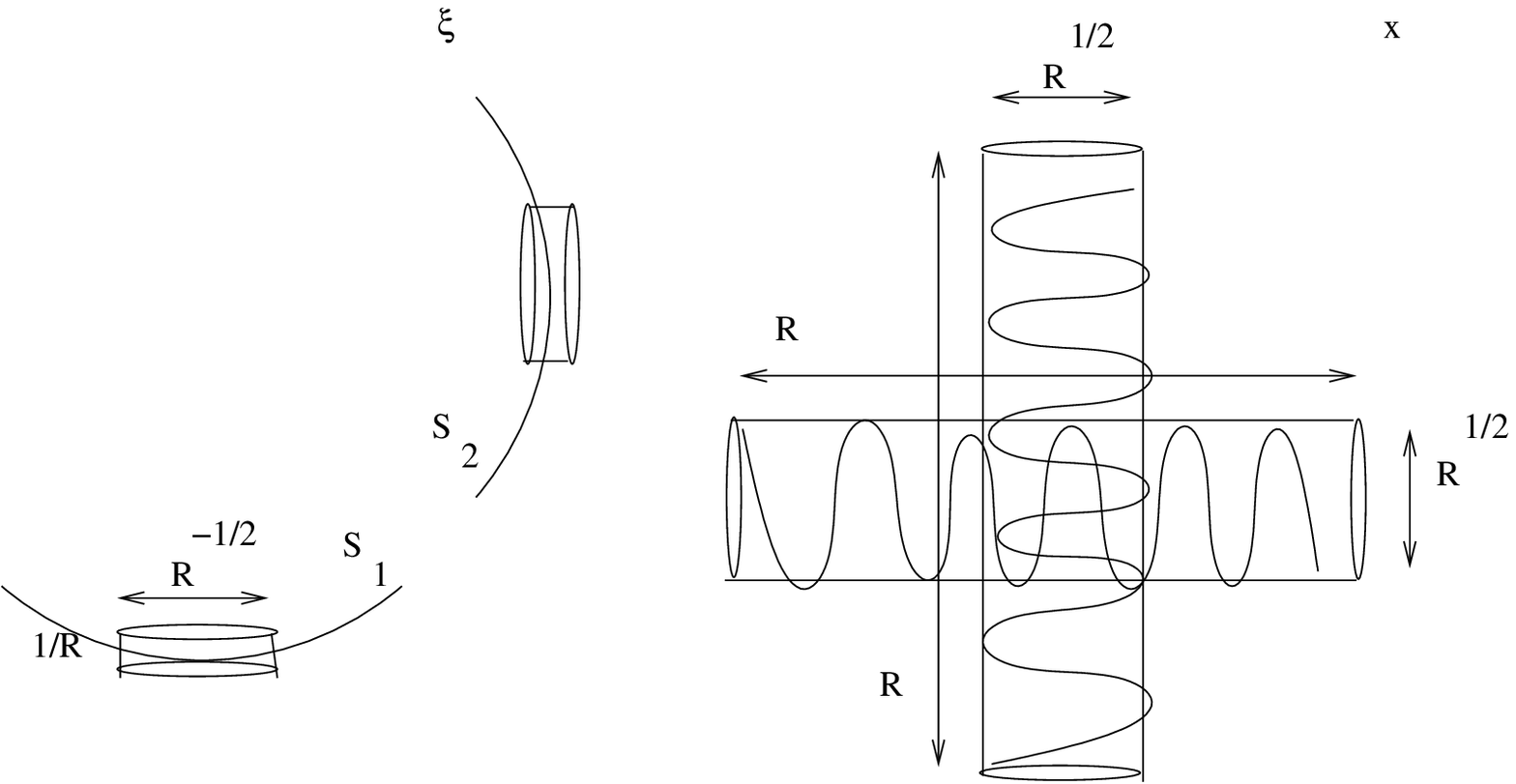}
 \caption{A naive generalization of the linear Knapp example to the bilinear setting.  This example is not very strong as the product of $(F_1 d\sigma_1)^\vee$ and $(F_2 d\sigma_2)^\vee$ will be much smaller in $L^p$ norms than is predicted by H\"older's inequality, because the supports of the two tubes depicted here are mostly disjoint.
}
\end{figure}

Now let $S_1$ and $S_2$ be two compact transverse subsets of the sphere, paraboloid, or cone.  Of course, any restriction theorem $R^*_S(q' \to p')$ for the sphere, paraboloid, or cone will imply a bilinear restriction theorem $R^*_{S_1,S_2}(q' \times q' \to p'/2)$ for the pair $S_1$, $S_2$.  However, the transversality allows us to prove more estimates in the bilinear setting than the linear one; we have already seen the estimate $R^*_{S_1,S_2}(2 \times 2 \to 2)$, whereas the linear restriction estimate $R^*_S(2 \to 4)$ is only true in three and higher dimensions.  One reason for this is that the Knapp example is much weaker in the transverse bilinear setting.  More precisely, the bilinear versions of the Knapp example will give substantially weaker necessary conditions for the bilinear restriction estimate $R^*_{S_1, S_2}(q' \times q' \to p'/2)$ than the Knapp condition $\frac{n+1}{p'} \leq \frac{n-1}{q}$ that one obtains for the linear restriction estimate $R^*_S(q' \to p')$.  Indeed, if one tries to directly mimic the Knapp construction, by setting $F_1$ to be the characteristic function of one spherical cap in $S_1$ of area $R^{-(n-1)}$, and $F_2$ to be another spherical cap in $S_2$ of area $R^{-(n-1)}$, then $\widehat{F_1 d\sigma_1}$ and $\widehat{F_2 d\sigma_2}$ are both functions of size $\sim R^{-(n-1)}$ on $R \times \R^2$ tubes $T_1$, $T_2$ of volume $\sim R^{n+1}$, as in the linear Knapp example; however, since $S_1$ and $S_2$ are transverse, $T_1$ and $T_2$ will be oriented at vastly different angles, and their intersection $T_1 \cap T_2$ will have much smaller volume, namely $\sim R^n$.  This will eventually lead to a weaker condition $\frac{n}{p'} \leq \frac{n-1}{q}$.  This is not the best known necessary condition; by modifying the bilinear Knapp example slightly one can obtain the conditions 
\begin{equation}\label{eq:bil-nec}
 p > \frac{2n}{n+1}; \quad \frac{n+2}{p'} + \frac{n}{q'} \leq n; \quad \frac{n+2}{p'} + \frac{n-2}{q'} \leq n-1;
\end{equation}
see \cite{tvv:bilinear}, \cite{damiano:null}, or Problem 2.8.  This is somewhat less stringent than the corresponding conditions
\begin{equation}\label{eq:knapp}
 p < \frac{2n}{n+1}; \quad \frac{n+1}{p'} + \frac{n-1}{q'} \leq n-1
\end{equation}
for the linear problem $R^*_S(q' \to p')$.  The \emph{bilinear restriction conjecture} asserts that the necessary conditions \eqref{eq:bil-nec} are in fact sharp.  This conjecture is still open except in dimension two, but recently it has been shown that (up to endpoints) it is equivalent to the usual restriction conjecture for the sphere, paraboloid, or cone \cite{wolff:cone}, \cite{tao:cone}, \cite{tao:parabola}; more in this in later lectures.

\begin{figure}[htbp] \centering
 \ \epsfbox{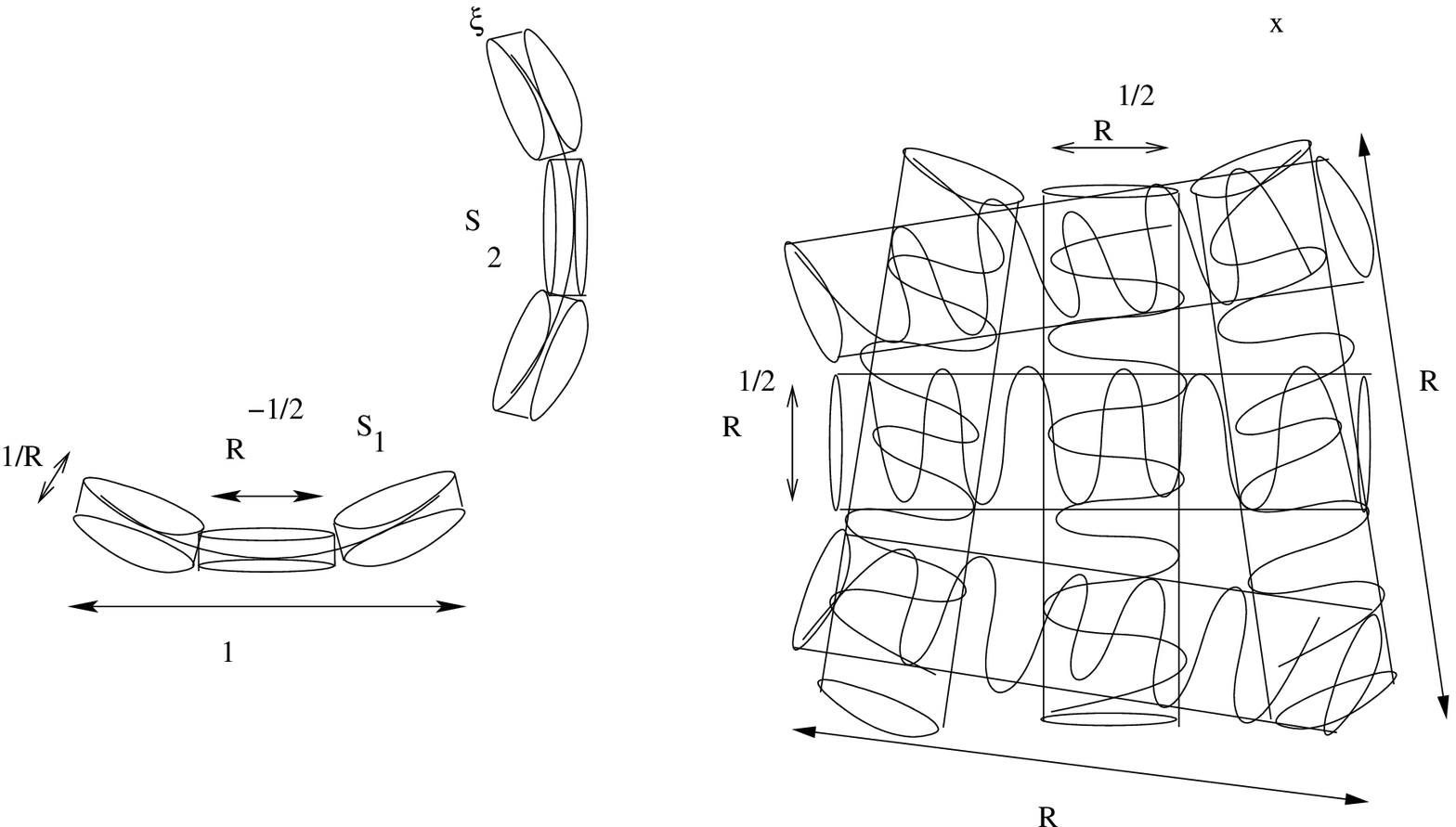}
 \caption{A better bilinear generalization of the linear Knapp example is the so-called ``stretched caps'' example, where several Knapp examples from both $S_1$ and $S_2$ are superimposed so that the supports of $(F_1 d\sigma_1)^\vee$ and $(F_2 d\sigma_2)^\vee$ match more closely.  Notice that both $(F_1 d\sigma_1)^\vee$ and $(F_2 d\sigma_2)^\vee$ can be viewed as a linear combination of oscillating tubes (``wave packets''); this is in fact a rather general fact, and will be exploited further in later sections.
}
\end{figure}

\begin{figure}[htbp] \centering
 \ \epsfbox{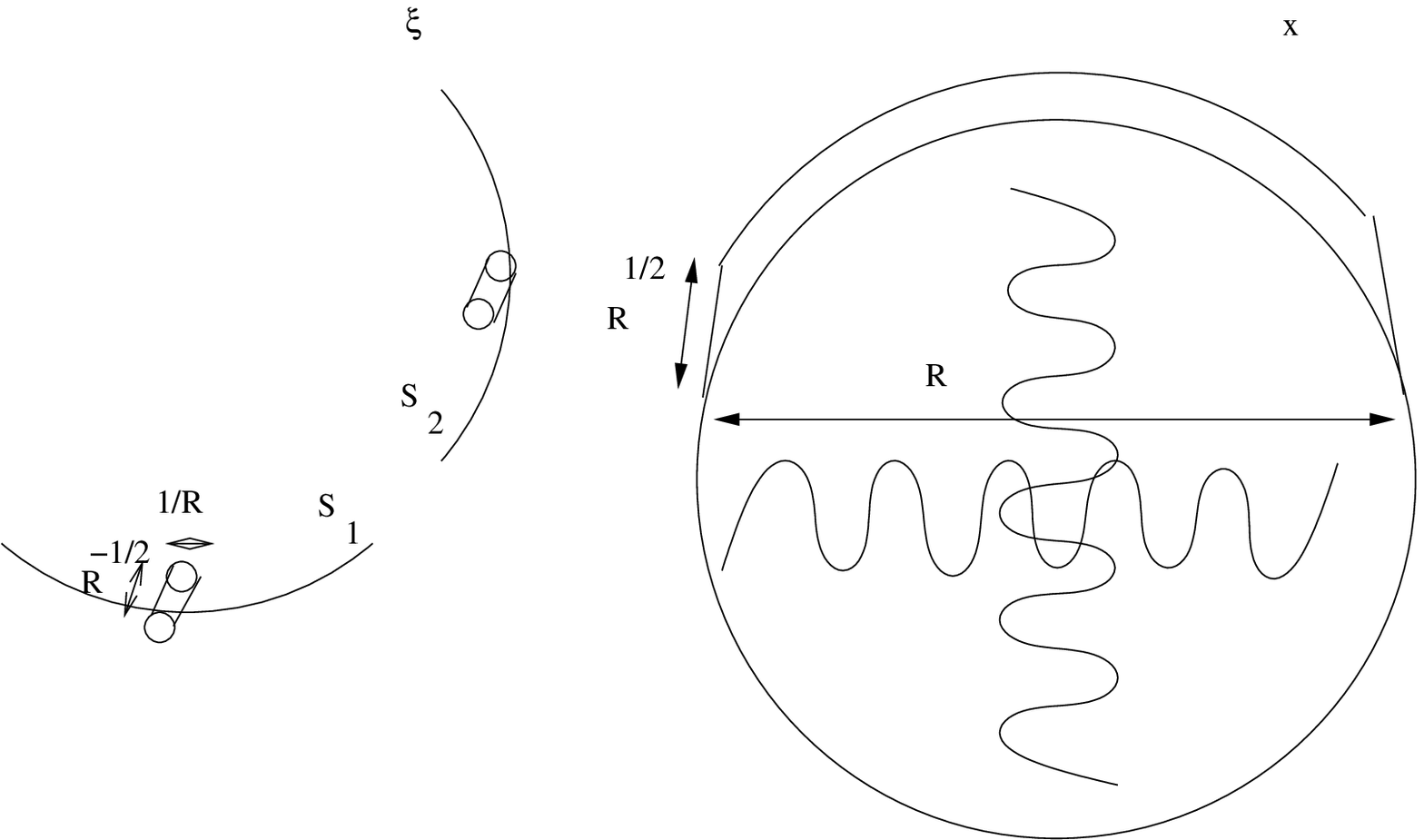}
 \caption{A variant of the ``stretched caps'' example is the ``squashed caps'' example, where $F_1$ and $F_2$ are now supported in thin parallel tubes with two dimensions comparable to $1/R$ and the rest comparable to $1/R^{1/2}$.  This places $(F_1 d\sigma_1)^\vee$ and $(F_2 d\sigma_2)^\vee$ both in the same disk (but with differing frequencies) with two dimensions comparable to $R$ and the rest comparable to $R^{1/2}$.  This example is superior to the squashed cap example when $q'$ is less than 2, as the $L^{q'}$ norms of $F_1$ and $F_2$ here are quite small.  One can partition this disk into tubes and thus view these examples also in the context of superpositions of wave packets.
}
\end{figure}

Up until now, we have viewed bilinear restriction estimates as being more complex generalizations of linear restriction estimates, which seems to offer no incentive to study the bilinear estimates until the linear ones are settled.  However, it turns out that one can use the bilinear estimates to go back and deduce new linear estimates, and indeed all the recent progress on the restriction problem has been obtained in this manner.  The key observation is that  one can perform a Whitney decomposition of the product manifold $S \times S$ around the diagonal $\Delta := \{ (\xi,\xi): \xi \in S\}$ so that $(S \times S) \backslash \Delta$ decomposes as the disjoint union of sets of the form $S_1 \times S_2$, where $S_1, S_2$ are disjoint subsets of $S$ whose separation is comparable to their diameter.  This allows one to obtain bilinear restriction estimates of the form $R^*_{S,S}(p' \times p' \to q'/2)$ (and hence $R^*_S(p' \to q')$) from estimates of the form $R^*_{S_1,S_2}(p' \times p' \to q'/2)$, using some rescaling and orthogonality estimates to sum up (and discarding the diagonal $\Delta$, which is of measure zero); see \cite{tvv:bilinear} for more details.  Of course one cannot hope to have an unconditional implication of the form $R^*_{S_1,S_2}(p' \times p' \to q'/2) \implies R^*_S(p' \to q')$, since the necessary conditions \eqref{eq:bil-nec} for the former are weaker than those \eqref{eq:knapp} for the latter; however, we can do the next best thing:

\begin{theorem}\label{thm:tvv}\cite{tvv:bilinear}  Let $p$, $q$ obey the necessary conditions \eqref{eq:knapp}, and suppose that $R^*_{S_1,S_2}(\tilde p' \times \tilde p' \to \tilde q')$ is true for all $(\tilde p, \tilde q)$ in an open neighborhood of $(p,q)$, and some pair $S_1, S_2$ of compact transverse subsets of the paraboloid.  Then $R^*_S(p' \to q')$ is true.
\end{theorem}

A similar result is true for the sphere, except that one must make $S_1$ and $S_2$ subsets of a certain parabolically rescaled version of the sphere; see \cite{tvv:bilinear} for more details.

The above theorem (and ones like it) allow one to pass freely back and forth between linear and (transverse) bilinear restriction estimates.  For instance, this theorem can be used to prove Zygmund's estimate (which asserts in particular that $R^*_S(4 \to 4+\eps)$ when $S$ is the unit circle) from the basic estimate $R^*_{S_1,S_2}(2 \times 2 \to 2)$ for transverse sets (this is basically the approach carried out in Problem 2.7).  Although the bilinear estimates appear more complicated, they are in fact easier to analyze because they consist purely of \emph{transverse interactions}, excluding the \emph{parallel interactions} which often cause the most trouble (cf. the Knapp example).

One can of course formulate local bilinear restriction estimates $R^*_{S_1,S_2}(p' \times p' \to q'; \alpha)$, which assert that
$$
\| (F_1 d\sigma_1)^\vee (F_2 d\sigma_2)^\vee \|_{L^{p'/2}(B(x_0,R))} \leq C_{p,q, S_1, S_2, \alpha} R^\alpha \| F_1 \|_{L^{q'}(S_1; d\sigma_1)} \| F_2 \|_{L^{q'}(S_2; d\sigma_2)}.
$$
One can of course reformulate these estimates using the uncertainty principle in a similar way to before, though some reformulations are not available because the notion of dualizing a bilinear estimate becomes difficult to use.  There are also bilinear ``epsilon-removal'' lemmas available; for instance, we have

\begin{theorem}\label{thm:eps}\cite{borg:cone}, \cite{tv:cone1}  Let $S_1$, $S_2$ be compact surfaces obeying some decay estimate 
$$ |(d\sigma_1)^\vee(x)|, |(d\sigma_2)^\vee(x)| \leq C/(1+|x|)^\rho$$
for some $\rho > 0$.  Suppose we have the local bilinear restriction estimate $R^*_{S_1,S_2}(2 \times 2 \to q, \eps)$ for all $\eps > 0$.  Then we have the global bilinear restriction estimate $R^*_{S_1,S_2}(2 \times 2 \to q+\eps)$ for all $\eps > 0$.
\end{theorem}

More quantitative versions of this estimate have been proven, see e.g. \cite{tv:cone1}, Lemma 2.4.  See also \cite{krt} for a more PDE-based approach to this epsilon-removal lemma.

The bilinear estimate $R^*_{S_1,S_2}(2 \times 2 \to 2)$ holds for all surfaces $S_1$, $S_2$ which are transverse.  If both $S_1$ and $S_2$ are flat, then this estimate is sharp; however one can improve this estimate slightly when $S_1$ and $S_2$ have some curvature.  For instance, if $S_1$ and $S_2$ are transverse subsets of the paraboloid in $\R^n$, then we have $R^*_{S_1,S_2}(p \times p \to 2)$ for all $p \geq \frac{4n}{3n-2}$, see \cite{tvv:bilinear}.  To see why we should gain over the $p=2$ estimate, consider the following.  Using Plancherel, we can rewrite $R^*_{S_1,S_2}(2 \times 2 \to 2)$ as
the bilinear convolution estimate
$$ \| (F_1 d\sigma_1) * (F_2 d\sigma_2) \|_{L^2(\R^n)}^2 \leq C \| F_1 \|_{L^2(S_1)}^2 \| F_2 \|_{L^2(S_2)}^2.$$
Let us suppose for the moment that we are in a model case, where the $F_j$, $j=1,2$ are characteristic functions, for the sets $\{ (\underline{\xi_j}, \frac{1}{2} |\underline{\xi_j}|^2): \underline{\xi_j} \in \Omega_j \}$ for some disjoint bounded open subsets $\Omega_1$, $\Omega_2$ of $\R^n$; the right-hand side is thus $C |\Omega_1| |\Omega_2|$.  Then by discarding some Jacobian factors (which are harmless due to the transversality), the left-hand side is essentially the volume of the $3n-1$-dimensional set
$$ \{ (\underline{\xi}_1, \underline{\xi}_2, \underline{\xi}_3, \underline{\xi}_4) \in \Omega_1 \times \Omega_2 \times \Omega_1 \times \Omega_2:
 \underline{\xi}_1+\underline{\xi}_2 = \underline{\xi}_3 + \underline{\xi}_4; \quad |\underline{\xi}_1|^2 + |\underline{\xi}_2|^2 = |\underline{\xi}_3|^2 + |\underline{\xi}_4|^2\}.$$
The two constraints imply that $\underline{\xi}_1$, $\underline{\xi}_2$ and $\underline{\xi}_3$, $\underline{\xi}_4$ form the opposing diagonals of a rectangle.  In particular, $\underline{\xi}_3$ lies on the hyperplane $\pi(\underline{\xi}_1, \underline{\xi}_4)$ containing $\underline{\xi}_1$ and orthogonal to $\underline{\xi}_4 - \underline{\xi}_1$, and then $\underline{\xi}_2$ can be recovered from the other three frequencies by the formula $\underline{\xi}_2 = \underline{\xi}_3 + \underline{\xi}_4 - \underline{\xi}_1$.  Thus, by Fubini's theorem, the volume of the above set is bounded above by
$$ \int_{\Omega_1} \int_{\Omega_2} | \Omega_2 \cap \pi(\underline{\xi}_1, \underline{\xi}_4) |\ d\underline{\xi}_1 d\underline{\xi}_4$$
(discarding the constraint $\xi_2 \in \Omega_2$).  Since $\Omega_2$ is bounded, we may make the very crude estimate
\begin{equation}\label{eq:crude}
 | \Omega_2 \cap \pi(\underline{\xi}_1, \underline{\xi}_4) | \leq C,
\end{equation}
from which the desired bound of $C |\Omega_1| |\Omega_2|$ follows.

\begin{figure}[htbp] \centering
 \ \epsfbox{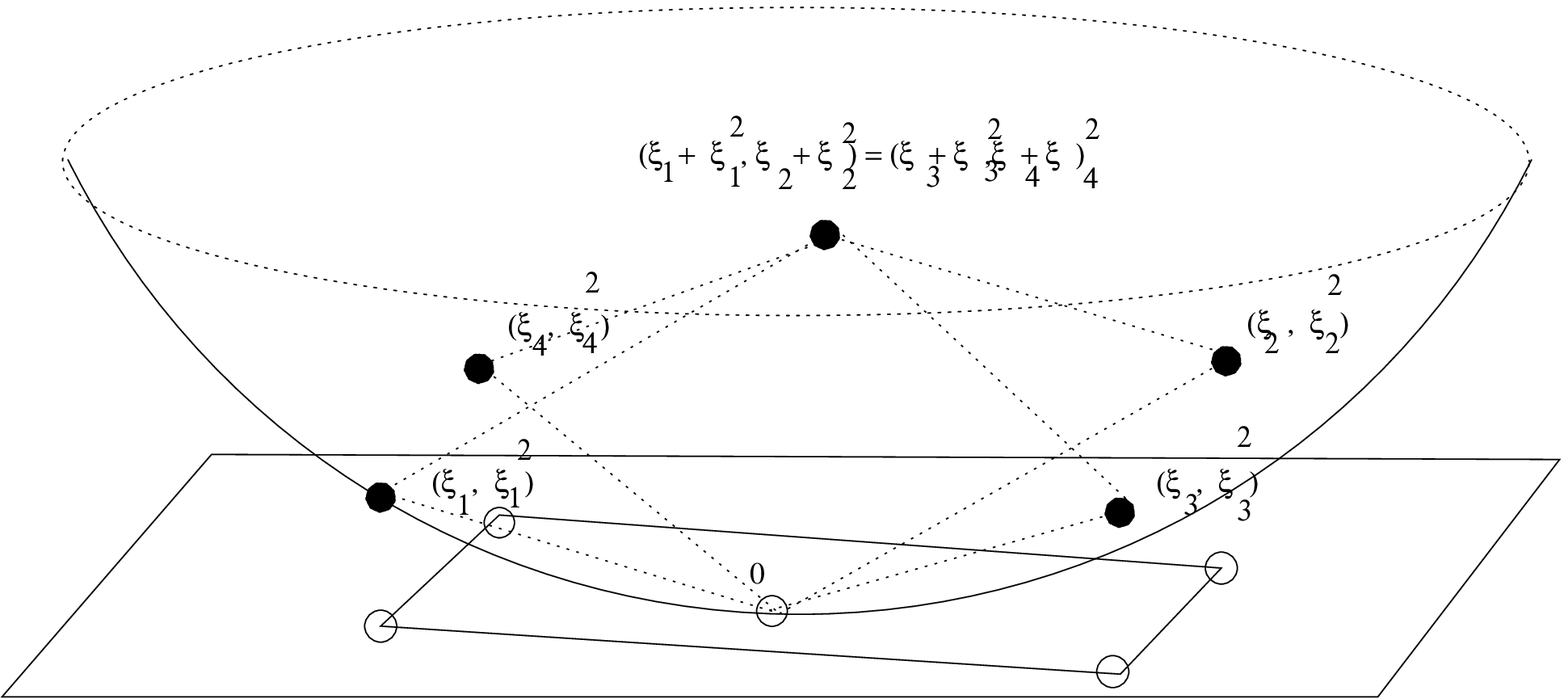}
 \caption{The $L^4$ geometry of the paraboloid.  In order for four points $(\underline{\xi_j}, |\xi_j|^2)$ to have equal sums, the co-ordinates $\underline{\xi_j}$ must be arranged in a rectangle.
}
\end{figure}

The estimate \eqref{eq:crude} can be improved (at least for generic values of $\underline{\xi}_1$, $\underline{\xi}_4$) when $\Omega_2$ has small measure, thanks the standard $L^p$ bounds for the Radon transform.  This is made rigorous in \cite{vargas:restrict}, \cite{vargas:2}, \cite{tvv:bilinear}, culminating in the above-mentioned bilinear restriction estimate $R^*_{S_1,S_2}(p \times p \to 2)$ for all $p \geq \frac{4n}{3n-2}$; this value of $p$ is best possible given that $q=2$, thanks to \eqref{eq:bil-nec}.  This issue of exploiting the possible gain over \eqref{eq:crude} also arises in some recent developments \cite{tao:parabola} in bilinear restriction theory, which we shall return to later.

The latter two conditions of \eqref{eq:bil-nec} meet when $q'=2$, when they assert that $R^*_{S_1, S_2}(2 \times 2 \to q)$ for all $q \geq \frac{2(n+2)}{n}$.  This was first conjectured by Machedon and Klainerman for both the paraboloid and cone.  Despite the original restriction conjecture remaining open, this conjecture has been completely solved for the cone and solved except for an endpoint for the paraboloid; see Figures 3, 4.  We shall discuss this recent progress in the next few lectures.

\begin{figure}[htbp] \centering
 \ \epsfbox{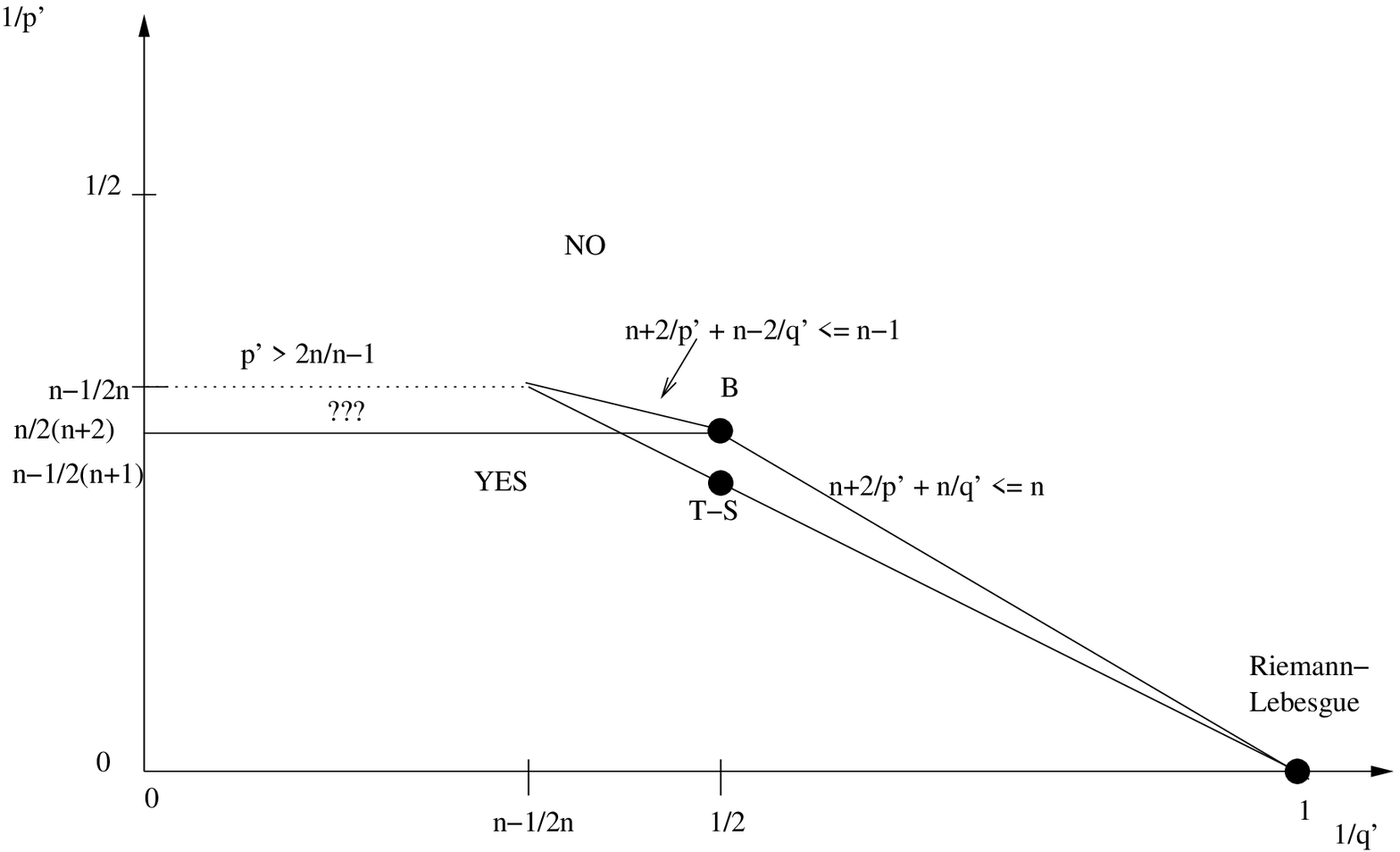}
 \caption{The restriction conjecture for transverse compact subsets of the sphere or paraboloid.  In addition to
the trapezoidal region in which linear restriction estimates are conjectured to hold, there is a larger pentagonal
region in which bilinear estimates are conjectured to hold.  The new vertex $B$ of this pentagon, where $(1/q',1/p') = (1/2, \frac{n}{2(n+2)})$, is still open, however we have restriction estimates at any exponent pair an epsilon below
this vertex \cite{tao:parabola}.  By H\"older this implies similar bilinear estimates to the left of $B$; using
the equivalence between linear and bilinear estimates in the trapezoidal region we can thus extend the region of known
restriction estimates to beyond the region given by the Tomas-Stein estimate.  (There were also several intermediate
results on this problem, see Figure \ref{fig:parab} and Figure \ref{fig:parab-bil}).
}
\end{figure}

\begin{figure}\label{fig:cone-bil}
\begin{tabular}{|l|l|l|} \hline
Dimension & Range of $q$ & \\

\hline
$n \geq 2$  & $q \geq 2$ & Plancherel + Cauchy-Schwarz\\ 
$n \geq 3$  & $q \geq \frac{n}{n-2}$ & Strichartz 1977 
\cite{strichartz:restrictionquadratic} \\
$n = 2$  & $q \geq 2 - \frac{13}{2408}$ & Bourgain 1995 \cite{borg:cone}\\ 

$n = 2$  & $q \geq 2 - \frac{8}{121}$ & Tao, Vargas 2000 \cite{tv:cone1}\\ 
$n \geq 2$ & $q > \frac{2(n+2)}{n}$ & Wolff, 2000 \cite{wolff:cone}\\
$n \geq 2$ & $q \geq \frac{2(n+2)}{n}$ & Tao, 2001 \cite{tao:cone} (best possible)\\
\hline

\end{tabular}
\caption{Known results on the bilinear restriction problem $R_{S_1 \times S_2}(2 \times 2 \to q)$, for transverse compact subsets of the cone. }
\end{figure}

\begin{figure}\label{fig:parab-bil}
\begin{tabular}{|l|l|l|} \hline
Dimension & Range of $q$ & \\
\hline

$n \geq 2$  & $q \geq 2$ & Plancherel + Cauchy-Schwarz\\ 
$n \geq 3$  & $q \geq \frac{n+1}{n-1}$ & Strichartz 1977 
\cite{strichartz:restrictionquadratic} \\
$n = 2$  & $q \geq 2 - \frac{5}{69}$ & Tao, Vargas, Vega 1998 \cite{tvv:bilinear}\\ 
$n = 2$  & $q \geq 2 - \frac{2}{17}$ & Tao, Vargas 2000 \cite{tv:cone1}\\ 

$n \geq 2$ & $q > \frac{2(n+2)}{n}$ & Tao, 2003 \cite{tao:parabola}\\
$n \geq 2$ & $q \geq \frac{2(n+2)}{n}$ & (conjectured)\\
\hline
\end{tabular}
\caption{Known results on the bilinear restriction problem $R_{S_1 \times S_2}(2 \times 2 \to q)$, for transverse compact subsets of the paraboloid. }
\end{figure}

\section{Problems for Lecture 2}

\begin{itemize}

\item {\bf Problem 2.1}.  Develop local analogues of the results in Problems 1.5 and 1.7.  (One could in fact develop local analogues of many more of the problems from that lecture, but this would get tedious after a while.  Note that for Problem 1.7, it may be convenient to use the formulation \eqref{eq:blur-adjoint}).

\item {\bf Problem 2.2}.(a)  Let $S$ be a smooth compact hypersurface, possibly with boundary, and let $d\sigma$ be surface measure on $S$.  Let $1 \leq p,q \leq \infty$.  Following the outline given in the notes, prove the equivalence of \eqref{eq:local-rest} and \eqref{eq:local-rest-blur}. (Hint:  To see how \eqref{eq:local-rest} implies \eqref{eq:local-rest-blur}, translate the surface $S$ by $O(1/R)$ and then average \eqref{eq:local-rest} over all such translations; to see the converse implication, introduce a bump function $\psi_{B(x_0,R)}$ concentrated near $B(x_0,2R)$ which equals 1 on $B(x_0,R)$, and exploit the reproducing formula $\hat f = \hat f * \hat \psi_{B(x_0,R)}$ to control $\hat f$ on $S$ in terms of $\hat f$ on neighborhoods such as $N_{1/R}(S)$, exploiting the fact that $\hat \psi_{B(x_0,R)}$ will decay rapidly away from $B(0,1/R)$.)

\item (b)*.  Now suppose in addition that we have the condition $q \leq p',p$.  Improve Problem 2.2(a) by showing that \eqref{eq:local-rest-blur} is in turn equivalent to the global version
\begin{equation}\label{eq:global-rest-blur}
  \| \hat f \|_{L^q(N_{1/R}(S))} \leq C_{p,q,S,\alpha} R^{\alpha-1/q} \| f \|_{L^p(\R^n)}
\end{equation}
holding for all $f \in L^p(\R^n)$.  (Hint: This is another manifestation of the uncertainty principle.  The implication of \eqref{eq:local-rest-blur} from \eqref{eq:global-rest-blur} is trivial.  For the converse implication,
use a smooth partition of unity to divide a global $f \in L^p(\R^n)$ into functions in $L^p(B(x_0,R))$ for various balls $B(x_0,R)$ and applying \eqref{eq:local-rest-blur} to each piece.  To sum, one will have use a reproducing formula as in (a) and to use an estimate of the form
\begin{equation}\label{hy}
 \| \sum_B F_B * \hat \psi_B \|_{L^q(\R^n)} \lesssim (\sum_B \| F_B \|_{L^q(\R^n)}^{\min(q,q')})^{1/\min(q,q')}
\end{equation}
for all $1 \leq q \leq \infty$, where $B$ ranges over a finitely overlapping collection of balls of radius $R$ and for each $B$, $\psi_B$ is a bump function adapted to $B$.  The estimate \eqref{hy} is a localized form of the Hausdorff-Young inequality and can be proven by complex interpolation starting from the extreme cases $q=1,2,\infty$.)

\item {\bf Problem 2.3}.  Again, let $S$ be a smooth compact hypersurface with boundary, and assume the ``Hausdorff-Young'' condition $q \leq p'$.  Show the equivalence of \eqref{eq:blur-adjoint} and \eqref{eq:discrete-adjoint}.  (Hint: This is yet another manifestation of the uncertainty principle.  To obtain \eqref{eq:discrete-adjoint} from \eqref{eq:blur-adjoint}, write the expression inside the norm of the left-hand side of \eqref{eq:discrete-adjoint} as $F$, and estimate $\| F \|_{L^{p'}(B(x_0,R))}$ by $\| F \psi_R \|_{L^{p'}(B(x_0,R))}$ for a suitable function $\psi_R$ whose Fourier transform is an approximation to the identity supported on the ball $B(0,1/R)$.  To prove the converse implication, cover $N_{1/R}(S)$ by translates of $\Lambda$ by $O(1/R)$, and thus write $G^\vee$ as an average of modulated versions of the expression in \eqref{eq:discrete-adjoint}, using a partition of unity if necessary.  Then use Minkowski's inequality.)

\item {\bf Problem 2.4.}  (a) Let $R \gg 1$, and let $\psi$ be a radial bump function adapted to the annular region $N_{1/R}(S_{sphere}) = \{ \xi \in \R^n: 1-\frac{1}{R} \leq |\xi| \leq 1 + \frac{1}{R} \}$; by this we mean that $\psi(x)$ depends only on the magnitude $r := |x|$ of $x$, is supported on the annulus $N_{1/R}(S_{sphere})$, and obeys the estimates
$$ \sup_{r > 0} |\partial_r^k \psi(r)| \leq C_{k} R^k$$
for all $k \geq 0$, where we have abused notation and written $\psi(x) = \psi(|x|) = \psi(r)$.  Show that $\psi^\vee(x) = O(R^{-(n-1)/2})$ for all $|x| \sim R$.  (Hint: decompose $\psi$ into bump functions adapted to disks of thickness $O(1/R)$ and radius $O(1/R^{1/2})$, obtain good estimates for the inverse Fourier transform of each piece, and then add up). 

\item (b) Conclude from this that $(d\sigma)^\vee(x) = O(\frac{1}{(1+|x|)^{(n-1)/2}})$ for all $x \in \R^n$, where $d\sigma$ is surface measure on the sphere.  (Hint: the case $|x| \leq 1$ is easy, so assume $|x| \sim R$ for some $R \geq 1$.  Then blur out $d\sigma$ by a radial approximation to the identity of width $1/R$ and use (a).  It is also possible to obtain this estimate, and indeed more precise asymptotics, via the method of stationary phase; however the point I wish to make here is that one can to a large extent duplicate the stationary phase computations by means of decompositions into Knapp examples, and heavy use of the uncertainty principle).

\item (c)*  Let $d\sigma$ be a smooth, compactly supported measure on $S_{cone}$ that avoids the origin. By modifying the above arguments, show that $(d\sigma)^\vee(x) = O(\frac{1}{(1+|x|)^{(n-1)/2}})$ for all $x \in \R^n$, where the implicit constants in the $O()$ notation depend of course on the exact choice of $d\sigma$.  (Hint: You will have to take the $O(1/R)$ neighborhood of the cone and divide it into ``slabs'' of length $O(1)$, angular width $O(1/R^{1/2})$, and thickness $O(1/R)$.)

\item {\bf Problem 2.5.}  (a) By squaring \eqref{eq:rpq}, show that the estimate $R_S(p \to 2)$ is equivalent to the estimate
$$ |\langle f * (d\sigma)^\vee, f \rangle| \lesssim \| f \|_{L^p(\R^n)}^2$$
holding for all $f \in L^p(\R^n)$.  Similarly, show that the estimate $R_S(p \to 2; \alpha)$ is equivalent to
\begin{equation}\label{eq:local-tts}
 |\langle f * (d\sigma)^\vee, f \rangle| \lesssim R^{2\alpha} \| f \|_{L^p(B(x_0,R))}^2
\end{equation}
holding for all $R \geq 1$, all balls $B(x_0,R)$ of radius $R$, and all functions $f \in L^p(B(x_0,R))$.

\item (b) Now suppose that $S$ is a compact smooth hypersurface with boundary, the measure $d\sigma$ obeys the decay estimate \eqref{eq:decay-ass}.  Prove that
\begin{equation}\label{eq:e-loc-glob}
 |\langle \chi_E * (d\sigma)^\vee, \chi_E \rangle| \lesssim R |E| + R^{-\rho} |E|^2
\end{equation}
for any measurable set $E$ and any $R \geq 1$.  (Hint: decompose $E$ into sets supported on balls of radius $R$; this decomposes the left-hand side of \eqref{eq:e-loc-glob} into a double sum.  Estimate the ``local'' part of this double sum using \eqref{eq:local-tts}, and the ``global'' part using the decay hypothesis.)  

\item (c) Optimize the estimate from (b) in $R$, and then reverse step (a), to obtain the estimate
$$ \| \hat{\chi_E}|_S \|_{L^2(S; d\sigma)} \lesssim \| \chi_E \|_{L^{2(\rho+1)/(\rho+2)}(\R^n)};$$
in other words, we have proven the Tomas-Stein estimate \eqref{ts} for characteristic functions.  Use Marcinkeiwicz interpolation to then conclude that $R_S(p \to 2)$ for all $1 \leq p < \frac{2(\rho+1)}{\rho+2}$.  (This is basically the Bourgain-Tomas version of the \emph{Tomas-Stein} argument; see \cite{tomas:restrict}, \cite{strichartz:restrictionquadratic}, \cite{borg:kakeya}, \cite{borg:stein} for more details.  One can obtain the endpoint $p = \frac{2(\rho+1)}{\rho+2}$ by using Stein's complex interpolation theorem; see \cite{stein:large}.  The proof of Theorem \ref{borg:global} begins by observing that the $|E|^2$ term on the right-hand side of \eqref{eq:e-loc-glob} can essentially be improved to $|E| \sup_{x_0 \in \R^n} |E \cap B(x_0,R)|$).

\item {\bf Problem 2.6.}   (a) Let $S$ be a $k$-dimensional compact hypersurface in $\R^n$, endowed with $k$-dimensional Hausdorff measure $d\sigma = d{\mathcal H^k}|_S$.  Show that a decay estimate of the form $(d\sigma)^\vee(\xi) = O( (1 + |\xi|)^{-\rho} )$ is only possible if $\rho \leq k/2$.  (Hint: Compute the $L^2$ norm of $d\sigma * \psi_{1/R}$, where $\psi_{1/R}$ is an approximation to the identity at scale $1/R$, and then use Plancherel to convert this to some information on $(d\sigma)^\vee$).

\item (b) Let $S$ be as in part (a).  Show that the restriction estimate $R^*_S(q' \to p')$ is only possible if $p' \geq \frac{2n}{k}$.

\item {\bf Problem 2.7.} (a) Let $S := \{ (x,y) \in \R_+ \times \R_+: x^2 + y^2 = 1 \}$ be a quarter-circle.  Let $R \geq 1$, and let $R^{-1/2} \leq \theta \lesssim 1$ be an angle.  Let $A$, $B$ be two arcs in $S$ of angle $\sim \theta$ and separation $\sim \theta$.  Show that
$$ \| \chi_{N_{1/R}(A)} * \chi_{N_{1/R}(B)} \|_\infty \lesssim R^{-2} \theta^{-1}.$$
Conclude (by an application of the Cauchy-Schwarz inequality) that
$$
\| G * H \|_{L^2(\R^n)} \lesssim R^{-1} \theta^{-1/2} \| G \|_{L^2(N_{1/R}(A))} \| H \|_{L^2(N_{1/R}(B))}$$
for all $L^2$ functions $G, H$ supported on $N_{1/R}(A)$ and $N_{1/R}(B)$ respectively.  (Hint: estimate $G * H$ pointwise by the geometric mean of $\chi_{N_{1/R}(A)} * \chi_{N_{1/R}(B)}$ and $|G|^2 * |H|^2$).  In particular, from H\"older's inequality, conclude that
\begin{equation}\label{eq:gh}
\| G * H \|_{L^2(\R^n)} \lesssim R^{-3/2} \| G \|_{L^4(N_{1/R}(A))} \| H \|_{L^4(N_{1/R}(B))}
\end{equation}

\item (b) For every function $G$ and $H$ defined on $N_{1/R}(S)$, and any $R^{-1/2} \leq \theta \lesssim 1$, define the partial convolution $G *_\theta H$ by

$$ G *_\theta H(x) = \int_{y+z = x; \angle y, z \sim \theta} G(y) H(z)\ dy.$$
Using \eqref{eq:gh}, show that
$$
\| G *_\theta G \|_{L^2(\R^n)} \lesssim R^{-3/2} \| G \|_{L^4(N_{1/R}(S))}^2
$$
for all $R^{-1/2} \leq \theta \lesssim 1$.   (Hint: Split $S$ into arcs $A$ of width $\theta$, and split $G$ accordingly; apply \eqref{eq:gh} to various pieces and then sum.  The key point here is that as one varies the arcs, the support of the corresponding portion of $G * _\theta G$ also varies, so that one has plenty of orthogonality). Conclude that \eqref{eq:gg} holds for $S$ equal to the quarter circle, and deduce the restriction theorem $R^*_S(4 \to 4, \eps)$ for every $\eps > 0$ for the quarter-circle, and hence for the circle.

\item {\bf Problem 2.8.}  Modify the Knapp example to obtain the necessary conditions \eqref{eq:bil-nec} for the bilinear restriction estimate $R^*_{S_1,S_2}(q' \times q' \to p'/2)$ for transverse subsets of the sphere.  (Hint: As remarked in the notes, a naive adaptation of the Knapp example will only yield the condition $\frac{n}{p'} \leq \frac{n-1}{q}$ this way.  To obtain the conditions \eqref{eq:bil-nec} you will have to 
align the Fourier supports of $\widehat{F_1 d\sigma_1}$ and $\widehat{F_2 d\sigma_2}$ better, either ``squashing'' the pair of caps or ``stretching'' them in some appropriate manner.)

\end{itemize}

\lecture{The wave packet decomposition}\label{sec:wave-packet}

We have discussed two of the tools in the modern theory of restriction estimates: the reduction to local estimates, and the reduction to bilinear estimates.  We now turn to a third key technique: the decomposition of $\widehat{fd\sigma}$ into \emph{wave packets}.

For sake of illustration, suppose we wish to prove the local restriction estimate $R_S(p \to 1; \alpha)$ where $S$ is the sphere \eqref{eq:sphere}; the exponent $1$ can of course be changed, but this does not significantly alter the argument sketched below (except that the estimates on certain coefficients $c_T$ will change).  We use the formulation \eqref{eq:blur-adjoint}, fixing $x_0 = 0$, thus we have to prove
$$\| G^\vee \|_{L^{p'}(B(0,R))} \lesssim R^{\alpha}$$
for all $R \geq 1$ and all functions $G$ in the unit ball of $L^\infty(N_{1/R}(S))$.  Henceforth we call a function ``bounded'' when it has an $L^\infty$ norm of $O(1)$.

Fix $R$ and $G$, and observe that the annular region $N_{1/R}(S)$ can be divided into $\sim R^{(n-1)/2}$ finitely overlapping disks $\kappa$ of width $\sim 1/\sqrt{R}$ and thickness $1/R$.  If $G$ is a function on $N_{1/R}(S)$, we can thus use a partition of unity to divide $G = \sum_\kappa G_\kappa$, where each $G_\kappa$ is a bounded function supported on one of these disks $G_\kappa$.  Our task is thus to show that
$$\| \sum_\kappa G_\kappa^\vee \|_{L^{p'}(B(0,R))} \lesssim R^{\alpha}.$$
The question then arises as to what $G_\kappa^\vee$ looks like.  We first consider some examples.  Suppose that the disk $G_\kappa$ is centered at a point $\omega_\kappa \in S^{n-1}$, which by the geometry of the sphere implies that $\omega_\kappa$ is also essentially the normal to the disk $\kappa$.  If $G_\kappa$ is a bump function adapted to $\kappa$, then by duality $G_\kappa^\vee$ would be concentrated on the $R \times R^{1/2}$ tube
$$ T_{\kappa,0} := \{ x \in B(0,R): \pi_{\omega_\kappa^\perp} x = O(R^{1/2}) \},$$
where $\pi_{\omega_\kappa^\perp}$ is the orthogonal projection onto the hyperplane $\omega_\kappa^\perp := \{ x \in \R^n: x \cdot \omega_\kappa = 0\}$.  Indeed, since $\kappa$ has volume roughly $R^{-(n+1)/2}$, we would expect $G_\kappa^\vee$ to equal a function $\psi_{T_{\kappa,0}}$ of the form
\begin{equation}\label{eq:wavepacket}
\psi_{T_{\kappa,0}}(x) = R^{-(n+1)/2} e^{2\pi i \omega_\kappa \cdot x} \phi_{T_{\kappa,0}},
\end{equation}
where $\phi_{T_{\kappa,0}}$ is a Schwartz function adapted to the tube $T_{\kappa,0}$ which has size $O(1)$ on this tube and is rapidly decreasing away from this tube.  We call the function $\psi_{T_{\kappa,0}}$ a \emph{wave packet} adapted to the tube $T_{\kappa,0}$; this object has already essentially come up in the discussion of the Knapp example.

What happens when $G_\kappa$ is not a bump function adapted to $\kappa$?  First suppose that $G_\kappa$ is a modulated bump function, more precisely suppose
$$ G_\kappa(\xi) = e^{-2\pi i x_0 \cdot \xi} \tilde G_\kappa(\xi)$$
where $\tilde G_\kappa$ is a bump function adapted to $\kappa$, and $x_0$ is an element of the hyperplane $\omega_\kappa^\perp$.  Then by the above discussion, $G_\kappa^\vee$ will be concentrated on the $R \times R^{1/2}$ tube
$$ T_{\kappa,x_0} := T_{\kappa,0} + x_0,$$
indeed we have
$$
G_\kappa^\vee(x) = \psi_{T_{\kappa,x_0}} := R^{-(n+1)/2} e^{2\pi i \omega_\kappa \cdot x} \phi_{T_{\kappa,x_0}}$$
for some Schwartz function $\psi_{T_{\kappa,x_0}}$ adapted to $T_{\kappa,x_0}$.
(One could also modulate $G_\kappa$ in the direction parallel to $\omega_\kappa$ instead of in the perpendicular directions, but this either has a negligible effect on the Fourier transform on the ball $B(0,R)$, or else makes the Fourier transform much smaller, depending on how much modulation is applied).

Thus one can make $G_\kappa^\vee$ resemble a wave packet $\psi_T$ for any tube $T$ oriented in the direction $\omega_\kappa$.  In the general situation, where $G_\kappa$ is a bounded function on $\kappa$, then one can perform a Fourier series decomposition in the directions perpendicular to $\omega_\kappa$ to essentially decompose $G_\kappa$ as an $l^2$-average of modulated bump functions.  (The behavior in the direction parallel to $\omega_\kappa$, which only extends for a distance $O(1/R)$ is essentially irrelevant, thanks to the localization of physical space to $B(0,R)$ and the uncertainty principle).  Thus we can write
\begin{equation}\label{eq:chop}
 G_\kappa^\vee = \sum_{T // \omega_\kappa} c_T \psi_T,
\end{equation}
where $T$ ranges over a finitely overlapping collection of $R \times \sqrt{R}$ tubes in $B(0,R)$ oriented in the direction $\omega_\kappa$, $\psi_T$ is a wave packet adapted to $T$, and $c_T$ is a collection of scalars with the $L^2$ normalization condition $\sum_{T // \omega_\kappa} |c_T|^2 \lesssim 1$; this can be thought of as a sort of windowed Fourier transform expansion for $G_\kappa$.
One can then expand the original Fourier transform $G^\vee$ as
$$ G^\vee = \sum_T c_T \psi_T$$
where $T$ now ranges over a separated\footnote{This means that any two tubes $T$, $T'$ in this collection either have directions differing by at least $1/R^{1/2}$, or are parallel and are separated spatially by at least $R^{1/2}$.} collection of tubes in $B(0,R)$, and $\omega_T$ denotes the direction of $T$.

This heuristic decomposition is an example of what is known as the \emph{wave packet decomposition} of $G^\vee$.  Versions of this decomposition in the context of the restriction problem (or the closely related Bochner-Riesz problem) first appeared in \cite{cordoba:covering}, \cite{cordoba:sieve}, \cite{feff:note}, \cite{feff:ball}, and was then later developed in \cite{borg:kakeya}, \cite{borg:stein}, \cite{vargas:restrict}, \cite{vargas:2}, \cite{tvv:bilinear}, \cite{tv:cone1}, \cite{wolff:cone}, \cite{tao:cone}; this method also can be applied to related problems such as local smoothing or Bochner-Riesz, see for instance \cite{wolff:distance}, \cite{wolff:smsub}.  A related, but slightly different, wave packet decomposition is also a standard tool in the analysis of Fourier integral operators (see e.g. \cite{sss}).  The wave packet decomposition reduces the study of restriction estimates to that of proving estimates on the linear superpositions of wave packets
\begin{equation}\label{eq:osc}
 \| \sum_T c_T \psi_T \|_{L^{p'}(B(0,R))}.
\end{equation}
Note that the wave packets $\psi_T$ have two main features; one at coarse scales $\gg \sqrt{R}$ and one at fine scales $\ll \sqrt{R}$.  At coarse scales, the wave packet is localized to a relatively thin tube of width $\sqrt{R}$ and length $R$.  At fine scales, the wave packet oscillates at a fixed frequency $\omega_T$.  Note that the coarse scale behavior and fine scale behavior are linked, because the direction of the tube at coarse scales is exactly the same as the frequency of the oscillation at fine scales.  The issue is then how to co-ordinate these two aspects - localization at coarse scales, and oscillations at fine scales - of wave packets in order to estimate \eqref{eq:osc} efficiently.

The first strategy for estimating these superpositions of wave packets is due to C\'ordoba \cite{cordoba:covering}, \cite{cordoba:sieve}, in which the idea is to estimate the oscillatory sum by the associated square function
\begin{equation}\label{eq:sq}
 \| (\sum_T |c_T \psi_T|^2)^{1/2} \|_{L^{p'}(B(0,R))}.
\end{equation}
The point of doing so is that all the fine-scale oscillation has been removed from this problem, leaving only the coarse scale localizations to tubes.  There is still of course the problem of estimating this non-oscillatory square function; this problem is essentially equivalent\footnote{Conversely, one must resolve the Kakeya conjectures in order to fully resolve the restriction problem, because one can use randomization arguments to show that any bound on \eqref{eq:osc} implies a comparable bound on \eqref{eq:sq}.  See e.g. \cite{bcss}, or Problem 3.2.} to the problem of estimating the \emph{Kakeya maximal function}, which is another important problem in harmonic analysis, but one which we will not discuss in detail here.  (See however \cite{wolff:kakeya}, \cite{Bo}, \cite{tao:elesc}, or Problems 3.2-3.3).

Now we discuss how to estimate the oscillatory sum \eqref{eq:osc} by the square function \eqref{eq:sq}.  When $p' = 2$, or when $n=2$ and $p' = 4$, one can bound the former by the latter by direct orthogonality (or bi-orthogonality) arguments (cf. Problem 2.6), however these arguments do not work for other values of $p'$.  Nevertheless, it was observed by Bourgain \cite{borg:kakeya}, \cite{borg:stein} that one can still obtain some control of \eqref{eq:osc} by \eqref{eq:sq} in these cases, but with a loss of some powers of $R$.  The idea is to break the ball $B(0,R)$ up into cubes $q$ of size $\sqrt{R}$.  On such ``fine-scale'' cubes, a wave packet $\psi_T = R^{-(n-1)/2} e^{2\pi i \omega_T \cdot x} \phi_T$ has essentially constant magnitude; to (over-)simplify the discussion, let us suppose that $\phi_T$ is equal to 1 on $q$ if $q \subset T$ and $\phi_T$ vanishes on $q$ otherwise.  Then the portion of \eqref{eq:osc} coming from $q$ is
$$ R^{-(n-1)/2} \| (\sum_{T: T \supset q} c_T e^{2\pi i \omega_T \cdot x} \|_{L^{p'}(B(0,R))}$$
while the corresponding portion of \eqref{eq:sq} is essentially
$$ R^{-(n-1)/2} R^{n/2p'} (\sum_{T: T \supset q} |c_T|^2)^{1/2}.$$
One can then control the former expression by the latter using discrete restriction estimates\footnote{It is intriguing that one uses local restriction estimates at scale $\sqrt{R}$, together with some Kakeya information, to obtain local restriction estimates at scale $R$.  This suggests a possible ``bootstrap'' approach where one could continually improve restriction estimates via iteration.  Some partial iteration methods to this effect can be found in \cite{borg:cone}, \cite{tvv:bilinear}, \cite{tv:cone1}; another example of this idea occurs in the induction-on-scales approach discussed in the next lecture.} of the type \eqref{eq:discrete-adjoint}, although the various powers of $R$ which accumulate when doing so do not necessarily all cancel, and so this method of estimation can cause some losses\footnote{It is conjectured that in any dimension $n \geq 2$, that one can estimate \eqref{eq:osc} by \eqref{eq:sq} in the endpoint case $p = 2n/(n+1)$, with at most an epsilon loss $R^\eps$; this, together with the so-called Kakeya maximal function conjecture \eqref{eq:kakeya}, would imply the restriction conjecture.  However, it is nowhere near solved at present, except when $n=2$, and is likely to be a harder problem than the restriction problem itself.}. 

By combining these observations with some non-trivial progress on the Kakeya maximal function conjecture, Bourgain \cite{borg:kakeya}, \cite{borg:stein} was able to obtain certain improvements to the Tomas-Stein estimate \eqref{ts} (see Figure \eqref{fig:parab}).  Further progress was made by Wolff \cite{wolff:kakeya}, who improved the Kakeya estimate used in Bourgain's argument.  By introducing bilinear (or $L^4$) methods to these arguments, further improvements were obtained in \cite{vargas:restrict}, \cite{tvv:bilinear}, \cite{borg:cone} \cite{tv:cone1}; one feature of these bilinear methods is that they could now be applied to the cone as well as the sphere or paraboloid.  

These methods, however, did not obtain sharp ranges of exponents, for a variety of technical reasons\footnote{The most obvious of these being that the Kakeya conjecture is still far from solved.  However, even if the Kakeya conjecture was completely resolved, there are still some remaining inefficiencies in the argument used to replace the oscillatory expression \eqref{eq:osc} by the square function \eqref{eq:sq}.}.  The next breakthrough was achieved by Wolff \cite{wolff:cone}, who solved (up to endpoints) the Machedon-Klainerman conjecture for cones, by employing one additional technique - that of induction on scales, which we discuss next.

\section{Problems for Lecture 3}

\begin{itemize}

\item {\bf Problem 3.1.} (a)  Let $F$ be a function supported on the cube $[1/3,2/3]^n$ in $\R^n$.  Show that there exists a decomposition of the inverse Fourier transform $F^\vee$ of the form
\begin{equation}\label{eq:unc}
 F^\vee = \sum_{k \in \Z^n} c_k \psi(x - k)
\end{equation}
where $\psi$ is a bump function on $\R^n$, and the $c_k$ are scalars such that
$$ \sum_{k \in \Z^n} |c_k|^2 \sim \|F\|_{L^2([1/3,2/3]^n)}.$$
(Hint: Expand $F$ as a Fourier series in $[0,1]^n$, and write $F = F \varphi$ for some bump function $\varphi$ supported on $[0,1]^n$ which equals 1 on $[1/3,2/3]^n$).  The decomposition \eqref{eq:unc} is, once again, another manifestation of the uncertainty principle - the Fourier transform of a function supported on what is essentially the unit cube, will itself be essentially constant on unit cubes.

\item (b) By translating and squashing $F$ by an appropriate linear transformation, use (a) to obtain a decomposition of the form \eqref{eq:chop} for Fourier transforms of bounded functions $G_\kappa$ supported on a $1/R \times 1/R^{1/2}$ spherical cap $\kappa$. 

\item {\bf Problem 3.2.}  (a) Suppose that we have a local restriction estimate $R^*_S(\infty \to 2p,\alpha)$ for the sphere for some $1 < p < \infty$.  Let $T_1, \ldots, T_K$ be any collection of $\sqrt{R} \times R$ tubes in $\R^n$, such that the directions $\omega_j \in S$ of the tubes $T_j$ are $R^{-1/2}$-separated (i.e. $|\omega_j - \omega_k| \geq R^{-1/2}$ whenever $j \neq k$).  Show that
$$ \| \sum_{j=1}^K \chi_{T_j} \|_{L^p(\R^n)} \lesssim 
R^{2\alpha + n-1}.$$
(Hint: start with \eqref{eq:blur-adjoint}, and apply this with $G = \sum_{j=1}^K \eps_j G_{\kappa_j}$, where $\eps_j = \pm 1$ are randomized signs and $G_{\kappa_j}$ are modulated Knapp examples adapted to a certain disk $\kappa_j$, designed so that $G_{\kappa_j}^\vee$ is large on $T_j$.  Then use Khinchin's inequality, see Appendix A.)

\item (b) Suppose that the restriction conjecture is true for the sphere (and in paritcular, that $R^*_S(\infty \to \frac{2n}{n-1},\eps)$ is true for any $\eps > 0$).  Deduce the estimate
\begin{equation}\label{eq:kakeya}
 \| \sum_{j=1}^K \chi_{T_j} \|_{L^{n/(n-1)}(\R^n)} \leq C_\eps \delta^{-\eps}
\end{equation}
for any $\eps > 0$ and any $0 < \delta \ll 1$, where $T_1, \ldots, T_K$ ranges over any collection of $\delta \times 1$ tubes whose directions $\omega_1, \ldots, \omega_K$ are $\delta$-separated.  (Hint: rescale the estimate obtained in (a)).  The estimate \eqref{eq:kakeya} is known as the \emph{Kakeya maximal function conjecture} and is an important unsolved problem in geometric combinatorics; it has been fully solved in dimension $n=2$, with partial progress in dimensions $n \geq 2$.  See e.g. \cite{wolff:kakeya}, \cite{Bo}, \cite{tao:notices}, \cite{tao:edinburgh} for more information on this and related problems.  The above result then asserts that the restriction conjecture implies the Kakeya conjecture; indeed, we strongly believe that one must first fully resolve the Kakeya conjecture before obtaining a full solution to the restriction conjecture\footnote{It is logically possible that one might be able to solve the restriction conjecture by other means and then deduce the Kakeya conjecture as a corollary, but this author doubts that this will be how the conjectures will be solved.  Also, in order to fully resolve the restriction conjecture, one will probably have to first prove not only the Kakeya conjecture, but various generalizations and improvements of that conjecture, possibly including, but not restricted to, bilinear variants, weighted variants, ``two-ends'' variants, or x-ray transform variants; several of these have already been employed to good effect on restriction problems in the literature, for instance the results discussed in the next lecture rely on bilinear, two-ends, direction-constrained variants of the Kakeya maximal function estimate.}.

\item {\bf Problem 3.3.}(a)  Assume that the Kakeya maximal function conjecture \eqref{eq:kakeya} is true.  Define a \emph{Besicovitch set} to be any subset $E$ of $\R^n$ such that $E$ contains a unit line segment in every direction (i.e. for every $\omega \in S_{sphere}$, there exists a unit line segment $l_\omega$ oriented in the direction $\omega$ and contained in $E$).  Show that for any $0 < \delta \ll 1$, the neighborhood $N_{\delta}(E)$ obeys the volume estimate
$$ |N_{\delta}(E)| \geq c_\eps \delta^\eps$$
for any $\eps > 0$.  (In particular, this implies that the \emph{Minkowski dimension}, also known as \emph{box counting dimension}, of $E$ is equal to $n$).  Hint: estimate $\langle \chi_{N_\delta(E)}, \sum_{j=1}^K \chi_{T_j} \rangle$ from above and below for a suitable collection $T_1, \ldots, T_K$ of tubes.

\item (b)* Let $E$ be a compact Besicovitch set, and assume that the map $\omega \to l_\omega$ from directions to unit line segments is Lebesgue measurable.  Again, we assume that the Kakeya maximal conjecture \eqref{eq:kakeya} is true.  Prove that $E$ has Hausdorff dimension equal to $n$.  Equivalently, for any $0 < d < n$ and every $A > 0$, show that there exists a $0 < \delta_0 \ll 1$ such that we have
$$ \sum_{j=1}^\infty r_j^d > A$$
for any covering $\{ B(x_j,r_j) \}_{j=1}^\infty$ of $E$ by balls $B(x_j,r_j)$ of radius $0 < r_j < \delta_0$.  (Hint: we may assume that all the $r_j$ are negative powers of 2.  For any $0 < \delta < \delta_0$ which is a negative power of 2, let $E_\delta := E \cap \bigcup_{j: r_j = \delta} B(x_j, r_j)$ and observe that the $E_\delta$ cover $E$.  Then let $0 < \eps \ll 1$ be a small number and use the pigeonhole principle to find a $\delta < \delta_0$ such that
$$ \int_{S_{sphere}} |l_\omega \cap E_\delta|\ d\omega \geq c_\eps \delta^\eps.$$
Then estimate $\langle \chi_{N_\delta(E_\delta)}, \sum_{j=1}^K \chi_{T_j} \rangle$ from above and below for a suitable collection $T_1, \ldots, T_K$ of tubes.)

\end{itemize}

\lecture{Induction on scales}\label{sec:induct}

The strategy to prove a local restriction estimate at a scale $R$ in the previous lecture can be summed up as follows: starting with a function $G^\vee$, decompose it into wave packets supported on $\sqrt{R} \times R$ tubes.  Designating scales greater than $\sqrt{R}$ as coarse, and scales less than $\sqrt{R}$ as fine, we use oscillatory estimates such as local restriction
estimates on fine scales, and Kakeya type estimates at coarse scales, in order to obtain the desired control on $G^\vee$.

This type of argument works particularly well when $G$ is a Knapp example supported on a disk of radius $R^{-1/2}$, so that $G^\vee$ is essentially a single wave packet.  However, it becomes inefficient when $G$ is a Knapp example spread out over a wider region, e.g. a cap-type region of radius $r^{-1/2}$ for some $1 \leq r \leq R$.  Then $G^\vee$ is concentrated on a much smaller set than a single wave packet - indeed, it is (somewhat) localized to a $r \times \sqrt{r}$ tube instead of an $R \times \sqrt{R}$ tube - but the wave packet decomposition requires that one decompose $G^\vee$ as the sum of much larger objects.  This is a rather inefficient decomposition, and one which leads to significant losses in the estimates.

The difficulty here is that the wave packet decomposition is chosen in advance, instead of being adapted to the particular function $G$ being investigated.  In particular, if it turns out that $G^\vee$ is concentrating in a much smaller region, say a ball $B(x_0,r)$, then one should replace the rather coarse $R \times \sqrt{R}$ wave packet decomposition by a finer one, in this case a $r \times \sqrt{r}$ decomposition.

Of course, the difficulty is that it would be incredibly complicated to actually try to construct such an adaptive wave packet decomposition, recursively passing from coarser scales to finer scales.  Fortunately, a way out of this complexity was discovered by Wolff \cite{wolff:cone} - which is to hide all this recursive complexity in an induction hypothesis, which we now refer to as an \emph{induction on scales} argument.  Using this new idea, Wolff was able to obtain a nearly-sharp bilinear restriction estimate for the cone, namely:

\begin{theorem}\label{wolff:thm} \cite{wolff:cone}
The bilinear restriction estimate $R^*_{S_1,S_2}(2 \times 2 \to q)$ is true for all transverse compact subsets $S_1$, $S_2$ of the cone $S_{cone}$ in $\R^n$, and all $q > \frac{n+2}{n}$.
\end{theorem}

The endpoint $q = \frac{n+2}{n}$ has since been obtained in \cite{tao:cone} by a refinement of the methods below.  This is sharp, see \eqref{eq:bil-nec}.

We now describe, rather informally, the idea of the arguments used to prove Theorem \ref{wolff:thm}; for a more rigorous presentation see \cite{wolff:cone}, \cite{tao:cone}, \cite{tao:non-endpoint}, \cite{tao:parabola}, \cite{krt}.  Suppose inductively that we already have some local estimate of the form $R^*_S(q' \to p'; \alpha)$; we will now try to use this estimate to prove a better estimate of the form $R^*_S(q' \to p'; \alpha - \eps)$ for some $\eps > 0$ depending on $\alpha$ (in what follows, the value of $\eps$ will vary from line to line).  Iterating this, we will eventually be able to obtain the estimate $R^*_S(q' \to p', \eps)$ for any $\eps > 0$, at which point we can use epsilon removal lemmas to obtain a global restriction estimate.

We still have to obtain the estimate $R^*(q' \to p'; \alpha - \eps)$ from the inductive hypothesis $R^*(q' \to p'; \alpha)$.  We first describe a somewhat oversimplified version of the main idea as follows.  We have to prove an estimate of the form
$$
\| (F d\sigma)^\vee \|_{L^{p'}(B(0,R))} \leq C_{p,q,S,\alpha} R^{\alpha-\eps} \| F \|_{L^{q'}(S)}
$$
for some $F$ on $S$ and $R \geq 1$, which we now fix.  Introduce the scale $r := R^{1-\eps}$, which is slightly smaller than $R$.  Then by the induction hypothesis $R^*(q' \to p'; \alpha)$ applied to scale $r$, we have
$$
\| (F d\sigma)^\vee \|_{L^{p'}(B(x_0,r))} \leq C_{p,q,S,\alpha} R^{\alpha-c\eps} \| F \|_{L^{q'}(S)}$$
for any ball $B(x_0,r)$.  Thus we can already prove the desired estimate on smaller balls $B(x_0,r)$.  More generally, we can prove
$$
\| (F d\sigma)^\vee \|_{L^{p'}(\bigcup_j B(x_j,r))} \leq C_{p,q,S,\alpha} R^{\alpha-c\eps} \| F \|_{L^{q'}(S)}$$
on any union $\bigcup_j B(x_j,r)$ of smaller balls, as long as the number of balls involved is not too large (e.g. at most $O((\log R)^C)$ for some absolute constant $C$).

As a rough first approximation, the idea of Wolff is to identify the ``bad'' balls $B(x_j,r)$ on which the function $(F d\sigma)^\vee$ ``concentrates''; the choice of these balls will of course depend on $F$.  These balls can be dealt with using the induction hypothesis, and it then remains to verify the restriction estimate on the exterior of these bad balls:
$$
\| (F d\sigma)^\vee \|_{L^{p'}(B(0,R) - \bigcup_j B(x_j,r))} \leq C_{p,q,S,\alpha} R^{\alpha-\eps} \| F \|_{L^{q'}(S)}.
$$

The above description of Wolff's argument was something of an oversimplification for two reasons; firstly, Wolff is working in the bilinear setting rather than the linear setting, and secondly the balls $B(x_j,r)$ turn out to depend not only on the original function $F$, but of the wave packet decomposition associated to $F$.  Let us ignore the first reason for the moment, and clarify the second.  On the ball $B(0,R)$, one can obtain a wave packet decomposition of the form
$$ (F d\sigma)^\vee(x) = \sum_T c_T \psi_T.$$
Because the argument of Wolff dealt with the cone, the wave packet decomposition here is slightly different from that discussed in the previous lecture, in two respects: firstly, the tubes $T$ are oriented on ``light rays'' normal to the cone $S$ instead of pointing in general directions, and secondly the internal structure of the wave packet $\psi_T$ is more interesting than just the product of a plane wave and a bump function, being decomposable into ``plates''.  We however will gloss over this technical issue.

For simplicity, let us suppose that the constants $c_T$ behave like a characteristic function; more precisely, there is some collection $\T$ of tubes such that $c_T = c$ for $T \in \T$ and $c_T = 0$ otherwise.  (The general case can be reduced to this case via a dyadic pigeonholing argument, which costs a relatively small factor of $\log R$).  Then we have
$$ (F d\sigma)^\vee(x) = c \sum_{T \in \T} \psi_T(x).$$
The idea now is to allow each wave packet $\psi_T$  to be able to ``exclude'' a single ball $B_T$ of the slightly smaller radius $r$.  In other words, one divides $(F d\sigma)^\vee$ into two pieces, a ``localized'' piece
$$ c \sum_{T \in \T} \psi_T(x) \chi_{B_T}(x)$$
and the ``global'' piece
$$ c \sum_{T \in \T} \psi_T(x) (1 - \chi_{B_T}(x)).$$
One then tries to control the localized piece using the induction hypothesis, and then handle the non-localized piece using the strategy of the previous lecture.

In the linear setting, this strategy does not quite work, because the localized pieces cannot be adequately controlled by the induction hypothesis.  However, in the bilinear setting, when one is trying to prove an estimate of the form
$$
\| (F_1 d\sigma_1)^\vee (F_2 d\sigma_2)^\vee \|_{L^{p'/2}(B(0,R))} \leq C_{p,q,S,\alpha} R^{\alpha-\eps} \| F_1 \|_{L^{q'}(S)} \| F_2 \|_{L^{q'}(S)}
$$
then one can decompose
$$ (F_j d\sigma_j)^\vee(x) = c_j \sum_{T_j \in \T_j} \psi_{T_j}(x)$$
for $j=1,2$, and allow each tube $T_j$ to exclude a single\footnote{Actually, in Wolff's argument there are $O((\log R)^C)$ such balls excluded, but this is a minor technical detail.} ball $B_{T_j}$ of radius $r$.  We can then split the bilinear expression
$$ \sum_{T_1 \in \T_1} \sum_{T_2 \in \T_2} \psi_{T_1} \psi_{T_2}$$
into a local piece
$$ \sum_{T_1 \in \T_1} \sum_{T_2 \in \T_2} \psi_{T_1} \psi_{T_2} \chi_{B_{T_1} \cap B_{T_2}}$$
(where \emph{both} tubes $T_1$ and $T_2$ are excluding $x$), and a global piece
$$ \sum_{T_1 \in \T_1} \sum_{T_2 \in \T_2}  \psi_{T_1} \psi_{T_2} (1 - \chi_{B_{T_1} \cap B_{T_2}}).$$
The local piece turns out to be easily controllable by the inductive hypothesis (the sum decouples into
non-interacting balls $B$ of radius $r$, and the claim follows by applying the hypothesis to each such ball separately
and summing), so it remains to control the global piece.

\begin{figure}[htbp] \centering
 \ \epsfbox{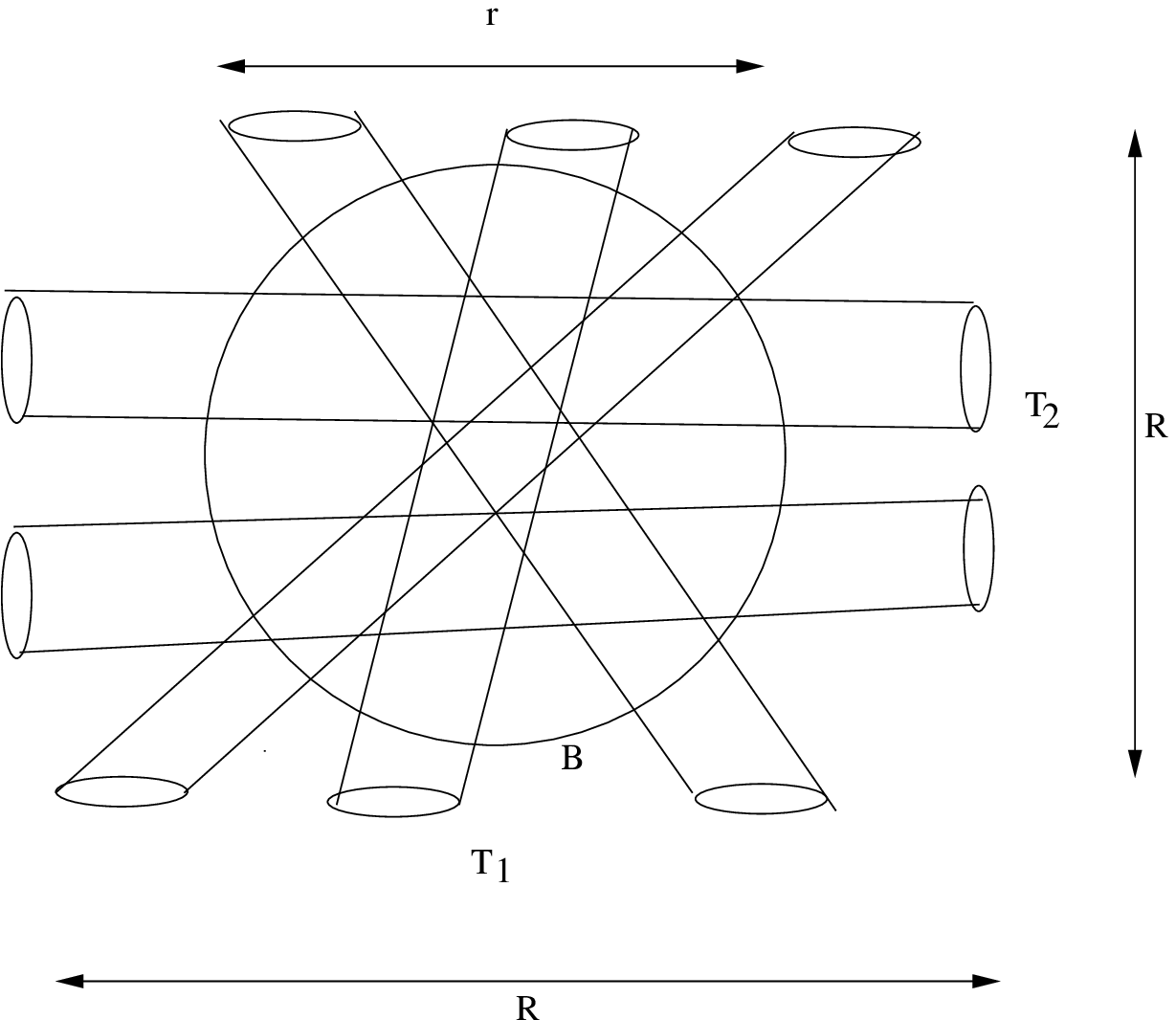}
 \caption{A very schematic depiction of the induction procedure.  The functions $(F_j d\sigma_j)^\vee$ are decomposed as the unions of wave packets supported on tubes.  Each tube $T$ is then associated with a ball $B_T$ of smaller radius $r < R$ where the tube has high interaction with other tubes; in the example given all the tubes would be associated
with a single ball $B_T = B$.  The local portion inside $B_T$ is estimated via the induction hypothesis, giving a constant of $r^\alpha = R^{\alpha - \eps}$ instead of $R^\alpha$.  For the global portion, the idea is to exploit the fact that there is significantly more disjointness between the tubes $T$ outside of the balls $B$.  Observe that in the squashed cap and stretched cap examples considered earlier, there is no advantage to be gained by excluding balls of radius $R$, because the incidences between tubes are well spread out.
}
\end{figure}

The key point is to prevent too many of the tubes $T_1$ and $T_2$ from interacting with each other.  This is done by selecting the balls $B_{T_1}$, $B_{T_2}$ strategically.  Roughly speaking, for each tube $T_1$, we choose $B_{T_1}$ to be the ball which contains as many intersections of the form $T_1 \cap T_2$ as possible; the ball $B_{T_2}$ is chosen similarly.  The effect of this choice is that any point $x$ which lies in a large number of tubes in $T_1$ and in $T_2$ simultaneously, is likely to be placed primarily in the local part of the bilinear expression, and not in the global part.

With this choice of the excluding balls $B_{T_1}$, $B_{T_2}$, Wolff was able to obtain satisfactory control on the number of times tubes $T_1$ from $\T_1$ would intersect tubes $T_2$ from $\T_2$.  The key geometric observation is as follows.  Suppose that many tubes $T_1$ in $\T_1$ were going through a common point $x_0$; since the tubes $T_1$ are constrained to be oriented along light rays, these tubes must then align on a ``light cone''.

Now consider a tube $T_2$ from $\T_2$; this tube is of course transverse to all the tubes $T_1$ considered above, and furthermore is transverse to the light cone that the tubes $T_1$ lie on.  It can either pass near $x_0$, or stay far away from $x_0$.  In the first case it turns out that the joint contribution of the tubes $T_1$ and $T_2$ will largely lie in the local part of the bilinear expression and thus be manageable.  In the second case we see from transversality that the tube $T_2$ can only intersect a small number of tubes $T_1$.

Thus if there is too much intersection among tubes $T_1$ in $\T_1$, then there will be fairly sparse intersection between those tubes $T_1$ and tubes $T_2$ in $\T_2$.  This geometric fact was exploited via combinatorial arguments in \cite{wolff:cone}, and when combined with some local $L^2$ arguments from \cite{mock:cone} to handle the fine scale oscillations, and the induction on scales argument, was able to obtain the near-optimal bilinear restriction theorem $R^*_{S_1, S_2}(2 \times 2 \to q)$ for $q > \frac{n+2}{n}$.  (The endpoint $q = \frac{n+2}{n}$ to the Machedon-Klainerman conjecture was then obtained in \cite{tao:cone} by refining the above argument).

\section{Adapting Wolff's argument to the paraboloid}\label{sec:parab}

The above argument of Wolff \cite{wolff:cone}, which yielded the optimal bilinear $L^2$ restriction theorem for the cone (Theorem \ref{wolff:thm}), relied on a key fact about the cone: all the tubes passing through a common point $x_0$, were restricted to lie on a hypersurface (specifically, the cone with vertex at $x_0$).  This property does not hold for the paraboloid, since in this setting the tubes can point in arbitrary directions.  Nevertheless, it is possible to recover this hypersurface property by exploiting a little more structure at fine scales, and more precisely by squeezing one ``dimension'' of gain out of \eqref{eq:crude}, thus obtaining the optimal bilinear $L^2$ restriction theorem for the paraboloid (and in fact also for the sphere, by a slight modification of the argument); this was achieved in \cite{tao:parabola}.  Specifically, we show

\begin{theorem}\label{tao:thm} \cite{tao:parabola}
The bilinear restriction estimate $R^*_{S_1,S_2}(2 \times 2 \to q)$ is true for all transverse (i.e. disjoint) compact subsets $S_1$, $S_2$ of the paraboloid $S_{parab}$ in $\R^n$, and all $q > \frac{n+2}{n}$.
\end{theorem}

We now give a sketch of this argument (or more precisely, how it differs
from Wolff's argument, on which it is based).
As in the last section, we can reduce matters to estimating a quantity such as
$$
\| \sum_{T_1 \in \T_1} \sum_{T_2 \in \T_2}  \psi_{T_1} \psi_{T_2} (1 - \chi_{B_{T_1} \cap B_{T_2}}) \|_{L^q(B(0,R))}$$
for some $1 < q < 2$.  It turns out in this case that one can obtain good bounds simply by interpolating between $L^1$ and $L^2$ bounds.  The $L^1$ bound is fairly trivial (using Cauchy-Schwarz to reduce to $L^2$ bounds on $\sum_{T_j \in \T_j} \psi_{T_j}$, which can be handled by orthogonality arguments), so we turn to problem of estimating the $L^2$ quantity:
$$
\| \sum_{T_1 \in \T_1} \sum_{T_2 \in \T_2}  \psi_{T_1} \psi_{T_2} (1 - \chi_{B_{T_1} \cap B_{T_2}}) \|_{L^2(B(0,R))}^2.$$
As is customary, we subdivide the large ball $B(0,R)$ into cubes $q$ of size $\sqrt{R}$.  The contribution of each cube $q$ is 
$$
\| \sum_{T_1 \in \T_1} \sum_{T_2 \in \T_2}  \psi_{T_1} \psi_{T_2} (1 - \chi_{B_{T_1} \cap B_{T_2}}) \|_{L^2(q)}^2.$$
Roughly speaking, we only need to consider pairs $T_1$, $T_2$ of tubes which pass through $q$ (because of the localization of $\psi_{T_1}$ and $\psi_{T_2}$, and such that $q$ is not contained in both $B_{T_1}$ and $B_{T_2}$.  For sake of argument, suppose that we only consider the terms where $q \not \subset B_{T_1}$.  Then we can rewrite the above expression as

$$
\| \sum_{T_1 \in \T'_1(q)} \psi_{T_1} \sum_{T_2 \in \T_2(q)} \psi_{T_2} \|_{L^2(q)}^2$$

where $\T_2(q)$ denotes all the tubes $T_2$ in $\T_2$ which intersect $q$, and $\T'_1(q)$ denotes all the tubes $T_1$ in $\T_1$ which intersect $q$ and for which $q \not \subset B_{T_1}$.  Note that the tubes in $\T'_1(q)$ must point in essentially different directions (since they all go through $q$, and are essentially distinct tubes), and similarly for $\T_2(q)$.

The function $\sum_{T_1 \in \T'_1(q)} \psi_{T_1}$ behaves roughly like the function $\frac{1}{R} (\chi_{\Omega'_1(q)} d\sigma_1)^\vee$, where $\Omega'_1(q)$ is the subset of the paraboloid whose unit normals lie within $1/R$ of the directions of one of the tubes in $\T'_1(q)$.  This can be seen by recalling the origin of these wave packets $\psi_{T_1}$, as Fourier transforms of functions on $S$ (or more precisely on $N_{1/R}(S)$; the discrepancy between the two explains the $\frac{1}{R}$ factor).  The function $\sum_{T_2 \in \T_2(q)} \psi_{T_2}$ is similarly comparable to the expression $\frac{1}{R} (\chi_{\Omega_2(q)} d\sigma_2)^\vee$ for a suitable set $\Omega_2(q)$.  Thus one is faced with an expression of the form
\begin{equation}\label{eq:22}
\| (\chi_{\Omega'_1(q)} d\sigma_1)^\vee (\chi_{\Omega_2(q)} d\sigma_2)^\vee \|_{L^2(\R^n)}^2,
\end{equation}
where we have discarded some powers of $R$ for sake of exposition, as well as the localization to $q$.  As observed previously, the restriction estimate $R^*(2 \times 2 \to 2)$ allows us to bound this quantity by something proportional to $|\Omega'_1(q)| |\Omega_2(q)|$.  Actually we can do a little better and refine this to, say, $|\Omega'_1(q)| |\Omega_2(q)|^2$; this comes from not discarding the constraint $\xi_2 \in \Omega_2$ in the argument immediately preceding \eqref{eq:crude}, and by exploiting the localization to $q$ more; we omit the details.  This is the type of bound used in Wolff's argument, in combination with the combinatorial arguments controlling the multiplicity of the tubes in $\T_1$ and $\T_2$ mentioned in the previous section, to obtain a sharp bilinear estimate in the case of the cone.  However, this bound is insufficient for the paraboloid case because of the failure of the tubes $T_1$ through a point to lie on a hypersurface.

Fortunately, this can be rectified by exploiting the gain inherent in \eqref{eq:crude}.  Indeed, by refusing to use \eqref{eq:crude} one can obtain a bound on \eqref{eq:22} which is proportional to
$$ |\Omega'_1(q)| |\Omega_2(q)| \sup_{\underline{\xi}_1, \underline{\xi}_4} |\Omega_2(q) \cap \pi(\underline{\xi}_1, \underline{\xi}_4)|.$$
This is similar to the bound of $|\Omega'_1(q)| |\Omega_2(q)|^2$ mentioned earlier, but is a little improved because one of the factors of $\Omega_2(q)$ is restricted to a hyperplane.  When one inserts this bound back into the coarse-scale combinatorial analysis of Wolff, this effectively allows us to restrict the tubes $T_2$ passing through a cube $q$ to be incident to a hyperplane.  This turns out to be a good substitute for the hypersurface localization property used in Wolff's argument, and is the key new ingredient which permits us to generalize the bilinear cone estimate to paraboloids (and by similar reasoning to other positively curved\footnote{It seems likely that this argument also extends to negatively curved surfaces, but this has not yet been checked rigorously.} surfaces, such as the sphere).  

One interesting feature of this argument is that it introduces a non-trivial correlation between the fine-scale analysis and the coarse-scale analysis; one may speculate that future developments on these problems will deal with the fine-scale and coarse-scale aspects of the restriction operator in a more unified manner. Certainly it seems that the interactions between multiple scales of Euclidean space will play a key role in the final resolution to the restriction problem.

\lecture{Connections with the Bochner-Riesz conjecture}

In our last two lectures, we discuss some connections between the restriction problem and some other open problems in harmonic analysis and PDE.  We have not attempted an exhaustive survey of the many interconnections, but hope to give a representative sample.  We have already mentioned the connection with the Kakeya conjecture, but we will not dwell more on this connection, despite it being a rich topic in itself, and refer the reader to other surveys \cite{wolff:kakeya}, \cite{Bo}, \cite{tao:elesc}, \cite{tao:notices}, \cite{tao:edinburgh} on this topic.

In the next lecture we shall discuss some connections between the restriction problem and PDE.  For this lecture we shall discuss instead the connection
with a number of oscillatory integral conjectures, focusing specifically on the \emph{Bochner-Riesz conjecture}.

Historically, the Bochner-Riesz conjecture arose from the classical problem \cite{bochner} of convergence of Fourier series and Fourier integrals in higher dimensions.  From Plancherel's theorem we know that for any $f \in L^2(\R^n)$, the Fourier transform $\hat f$ is also in $L^2(\R^n)$, and we have the Fourier inversion formula
$$ f(x) = \int_{\R^n} \hat f(\xi) e^{2\pi i x \cdot \xi}\ d\xi$$
holds in the following sense: if $K_j$ is any increasing set of compact sets which converges to $\R^n$, then the partial Fourier integrals
$$ \int_{K_j} \hat f(\xi) e^{2\pi i x \cdot \xi}\ d\xi $$
converge in $L^2(\R^n)$ norm to $f$ as $j \to \infty$.  In particular, the integrals
\begin{equation}\label{eq:disk}
 \int_{B(0,R)} \hat f(\xi) e^{2\pi i x \cdot \xi}\ d\xi 
\end{equation}
converge in $L^2(\R^n)$ norm to $f$ as $R \to \infty$.

We can rewrite this fact in the language of Fourier multipliers.  For any bounded function $m \in L^\infty(\R^n)$, define the Fourier multiplier $m(D): L^2(\R^n) \to L^2(\R^n)$ by the formula
$$ \widehat{m(D) f}(\xi) := m(\xi) \hat f(\xi).$$
Then the expression \eqref{eq:disk} is just $\chi_{B(0,R)}(D) f(x)$.  The operator $\chi_{B(0,R)}(D)$ is known as the \emph{ball multiplier} (or \emph{disk multiplier} in two dimensions) at frequency $R$.  Plancherel's theorem then asserts that for any $f \in L^2(\R^n)$, the functions $\chi_{B(0,R)}(D) f$ converge to $f$ in $L^2$ norm\footnote{It is an extremely interesting and difficult open problem as to whether one also has almost everywhere pointwise convergence of these means.  This is only known in one-dimension, in a celebrated theorem of Carleson \cite{carleson}.  In higher dimensions this result would be an endpoint version of the \emph{maximal Bochner-Riesz} conjecture, but it seems that this problem requires substantially new techniques beyond what is used for Bochner-Riesz problems, and in particular some careful and sophisticated time-frequency decompositions in higher dimensions.} 
as $R \to \infty$.

One may then ask whether the same result is true in $L^p(\R^n)$ for other values of $p$.  By the uniform boundedness principle, this is equivalent to the ball multipliers $\chi_{B(0,R)}(D)$ being bounded on $L^p(\R^n)$ uniformly in $R$.  But a simple scaling argument shows that the $(L^p, L^p)$ operator norm of $\chi_{B(0,R)}(D)$ is independent of $R$, so it suffices to show that the uniform ball multiplier $\chi_{B(0,1)}(D)$ is bounded on $L^p$.

In one dimension, $\chi_{B(0,1)}(D)$ is a linear combination of modulated Hilbert transforms, and hence  are bounded in $L^p$ for every $1 < p < \infty$ (as one can show by many methods, e.g. by Calder\'on-Zygmund theory).  However, in 1971 Fefferman \cite{feff:ball} proved the remarkable result that the ball multipliers are unbounded in $L^p(\R^n)$ for every $n \geq 2$ and $p \neq 2$.

The proof of this fact proceeds by a variant of the wave packet decomposition discussed in earlier lectures, and we sketch the proof heuristically as follows.  
By duality it suffices to consider the case $p>2$.  Let $R \gg 1$, and let $T$ be an $\sqrt{R} \times \R$ in $\R^n$, and oriented in some direction $\omega_T$.  Let $\psi_T$ be a bump function adapted to  the tube $T$, and let $\tilde T$ be a shift of $T$ by $2R$ units in the $\omega_T$ direction.  Then a computation shows that
$$ |\chi_{B(0,1)}(D) ( e^{2\pi i \omega_T \cdot x} \psi_T(x) )| \approx 1$$
for all $x \in \tilde T$.  Indeed, the Fourier transform of $ e^{2\pi i \omega_T \cdot x} \psi_T(x)$ is concentrated in a $\frac{1}{\sqrt{R}} \times \frac{1}{R}$ cap $\kappa_T$ in $N_{1/R}(S_{sphere})$ centered at $\omega_T$ and with unit normal $\omega_T$.  Applying the multiplier $\chi_{B(0,1)}(D)$ effectively bisects this cap $\kappa_T$ in half, which by the uncertainty principle will inevitably stretch the function $e^{2\pi i \omega_T \cdot x} \psi_T(x)$ in the longitudinal direction $\omega_T$.

We now apply considerations similar to that of Problem 3.2.  Let $T_1, \ldots, T_K$ be some collection of $\sqrt{R} \times R$ tubes, oriented in a $\frac{1}{R}$-separated set of directions, and $f$ be the function
$$ f(x) := \sum_{j=1}^K \eps_j e^{2\pi i \omega_{T_j} \cdot x} \psi_{T_j}(x)$$
where $\eps_j = \pm$ are randomized signs.  Then on the average, the $L^p$ norm of $f$ is roughly
$$ \| f \|_p \sim \| (\sum_{j=1}^K \chi_{T_j})^{1/2} \|_p$$
thanks to Khinchin's inequality, while on the average, the $L^p$ norm of $\chi_{B(0,1)}(D) f$ is bounded below by
$$ \| \chi_{B(0,1)}(D) f \|_p \gtrsim \| (\sum_{j=1}^K \chi_{\tilde T_j})^{1/2} \|_p.$$
Thus, to disprove the $L^p$ boundedness of $\chi_{B(0,1)}(D)$, it suffices to find, for each $A$, an arrangement of direction-separated tubes $T_1, \ldots, T_K$ such that
$$\| (\sum_{j=1}^K \chi_{\tilde T_j})^{1/2} \|_p \gtrsim A \| (\sum_{j=1}^K \chi_{T_j})^{1/2} \|_p $$
or equivalently
$$\| \sum_{j=1}^K \chi_{\tilde T_j} \|_{p/2} \gtrsim A^2 \| \sum_{j=1}^K \chi_{T_j} \|_{p/2}. $$
This cannot happen when $p=2$ (because then one can pull the sum out of the norm and observe that $\tilde T_j$ and $T_j$ have essentially the same volume), but when $p > 2$ it is possible to arrange the tubes $T_j$ so that they are disjoint, but that their shifts $\tilde T_j$ have high overlap, which causes the $p/2$ norm to be higher.  (This construction is detailed in \cite{feff:ball} or \cite{stein:large}, and is based on a variant of Besicovitch's construction of a Besicovitch set (see Problem 3.3(b)) with measure zero).

The failure of $L^p$ convergence of the ball multipliers can be worked around by replacing the ball multipliers with a smoother family of multipliers, known as the \emph{Bochner-Riesz multipliers}
$$ S^\delta_R := (1 - |\frac{D}{R}|^2)_+^\delta$$
where $x_+ := \max(x,0)$ denotes the positive part of $x$.  When $\delta = 0$, the multiplier $S^\delta_R$ is the same as $\chi_{B(0,R)}(D)$, but for $\delta > 0$ the multiplier is a little bit smoother at the boundary of $B(0,R)$ (for instance, one can easily represent $S^\delta_R$ as the average (modulo a tractable bump function error) of ball multipliers $\chi_{B(0,R')}(D)$ where $R'$ varies close to $R$).  For $\delta < 0$ the multiplier has an unbounded symbol, which means it cannot map $L^p$ to itself for any $p$ (although it is still possible to have $L^p$ to $L^q$ estimates, see \cite{bak}).
The question then arises: for which $p$ and which $\delta > 0$ do the means $S^\delta_R f$ converge to $f$ in $L^p$ norm as $R \to \infty$, for all $f \in L^p(\R^n)$?  By the same considerations as before, this is equivalent to having the estimate
\begin{equation}\label{eq:dp}
 \| S^\delta_1 f \|_{L^p(\R^n)} \lesssim \| f \|_{L^p(\R^n)}.
\end{equation}
This problem is known as the \emph{Bochner-Riesz problem}.  When $\delta$ is very large, this estimate is fairly easy, since the symbol of $S^\delta_1$ then becomes very smooth and one can apply for instance the H\"ormander multiplier theorem.  In fact, one can think of the Bochner-Riesz problem as an attempt to probe the sharp limits of the H\"ormander multiplier theorem, for the model case of multipliers singular on the sphere.

We have already observed the necessary condition $\delta > 0$ (except for the exceptional case $\delta = 0, p = 2$).  Another necessary condition, discovered by Herz, arises from considering the convolution kernel $K_\delta := ((1 - |\xi|^2)_+^\delta)^\vee$ of $S^\delta_1$.  Standard stationary phase asymptotics (see e.g. \cite{stein:large}; one can also proceed by dividing the symbol of $S^\delta_1$ into a large number of Knapp examples) give the asymptotics
$$ K_\delta(x) \sim e^{\pm 2\pi i |x|} / |x|^{(n+1)/2 + \delta}$$
which is only in $L^p(\R^n)$ when
$$ n(\frac{1}{2} - \frac{1}{p}) - \frac{1}{2} < \delta.$$
Applying $S^\delta_1$ to a bump function, we thus see that the above condition is a necessary condition in order for $S^\delta_1$ to be bounded on $L^p(\R^n)$.  
By duality, we thus obtain the necessary condition
$$ n|\frac{1}{2} - \frac{1}{p}| - \frac{1}{2} < \delta.$$
The \emph{Bochner-Riesz conjecture} then asserts that this condition, in addition to the condition $\delta > 0$ mentioned previously, are not only necessary, but also sufficient for $S^\delta_1$ to be bounded on $L^p(\R^n)$.  This conjecture has been verified in two dimensions but remains open (despite much partial progress) in higher dimensions.

One can simplify the estimate \eqref{eq:dp} using dyadic decomposition.  Observe that we may decompose
$$ (1 - |\xi|^2)^\delta_+ = m_0(\xi) + \sum_{j=1}^\infty 2^{-\delta j} m_j(\xi)$$
where $m_0$ is a bump function adapted to the unit ball, and $m_j(\xi)$ is a bump function adapted to the annular region $\{ \xi: |\xi| = 1 + O(2^{-j}) \}$.  (One can impose more conditions on $m_j$, which can be useful for the endpoint theory, see \cite{tao:weak}; however for this discussion the exact form of $m_j$ is irrelevant).  We can thus decompose
$$ S^\delta_1 = m_0(D) + \sum_{j=1}^\infty 2^{-\delta j} m_j(D).$$
The operator $m_0(D)$ is harmless (it is convolution with $m_0^\vee$, which is a Schwartz function).  Thus, giving up an epsilon in the $\delta$ parameter\footnote{The Bochner-Riesz conjecture is still far from completely solved, in that the best values of $\delta$ we can obtain are often quite far from the optimal ones conjectured.  Because of this, one is usually quite blas\'e about conceding epsilons in the $\delta$ parameter - which is similar to the $\alpha$ parameter in local restriction estimates.  In any event there are epsilon-removal lemmas which allow one to trade epsilon losses in the $\delta$ parameter for those in the $p$ parameter, see \cite{tao:weak2}.} if necessary, the Bochner-Riesz estimate \eqref{eq:dp} is thus equivalent to the estimate
\begin{equation}\label{eq:dp-local}
 \| m_j(D) f \|_{L^p(\R^n)} \lesssim 2^{\delta j} \| f \|_{L^p(\R^n)} \hbox{ for all } j \geq 1.
\end{equation}
This estimate should be compared with the restriction estimate \eqref{eq:global-rest-blur}.  In \eqref{eq:global-rest-blur}, one started with a function on the annulus $N_{1/R}(S_{sphere})$, and then measured its inverse Fourier transform in $L^p(\R^n)$.  Here, we are starting with a function in $L^p(\R^n)$, taking its Fourier transform, restricting that smoothly to the annulus $N_\delta(S_{sphere})$, and then taking inverse Fourier transforms again and measuring that in $L^p(\R^n)$.  Thus the estimate \eqref{eq:dp-local} differs from \eqref{eq:global-rest-blur} mainly by a Fourier transform.  This connection can be exploited efficiently when there is an $L^2$ norm lying around somewhere; for instance in \cite{feff:note} Fefferman observed that whenever one has a restriction theorem $R^*_S(2 \to p)$ for the sphere, this implies the Bochner-Riesz conjecture is true for that exponent $p$; this follows from the above connection, Plancherel, and a judicious application of the uncertainty principle.  Thus for instance, the Tomas-Stein restriction theorem \eqref{ts} implies the Bochner-Riesz conjecture for $p \geq \frac{2(n+1)}{n-1}$ (or $p \leq \frac{2(n+1)}{n+3}$, by duality).  Various improvements to this result where then developed by Bourgain \cite{borg:kakeya}, \cite{borg:stein} and others (see \cite{vargas:restrict}, \cite{vargas:2}, \cite{wolff:kakeya}, \cite{tvv:bilinear}, \cite{tv:cone1}).  More recently, Lee \cite{lee:2} has established the bilinear analogue of Fefferman's argument in that whenever the bilinear restriction theorem $R^*_{S_1,S_2}(2 \times 2 \to p/2)$ holds for transverse subsets $S_1, S_2$ of the sphere (or other positively curved hypersurfaces), then the Bochner-Riesz conjecture is true for that exponent $p$.    This (combined with Theorem \ref{tao:thm}) has led to the best known result on the Bochner-Riesz problem in three and higher dimensions, namely that the conjecture is true when $p \geq \frac{2(n+2)}{n}$ or $p \leq \frac{2(n+2)}{n+4}$.  (For intermediate values of $p$, the best result to date can be obtained by interpolating the above result with the $p=2$ theory).  

This is not the only connection between the restriction and Bochner-Riesz conjecture.  For instance, the Bochner-Riesz conjecture is known to imply the restriction conjecture (which is in some sense the scattering limit of the Bochner-Riesz conjecture); see \cite{tao:boch-rest}.  This implication can be reversed if the sphere is replaced by the paraboloid; see \cite{carbery:parabola}.

The Bochner-Riesz problem itself comes with many variants, including a weak-type endpoint (see e.g. \cite{tao:weak2} and the references therein), Hardy space versions (see e.g. \cite{stein:hp}, a maximal version (see e.g. \cite{christ:weight}, \cite{carbery:maximal-bochner}, \cite{tao:weak2}, \cite{tao:bochmax}), weighted versions (see e.g. \cite{stein:problem}, \cite{carbery:weight}), and versions for surfaces other than the sphere (see e.g. \cite{mock:cone}, \cite{borg:cone}, \cite{tv:cone2} for a discussion of the Bochner-Riesz problem for the cone).  The behavior at the endpoint $p=1, \delta = (n-1)/2$ has also been extensively studied, see e.g. \cite{liflyand}.  The Bochner-Riesz multipliers are also closely related to the fundamental solution operators for the Helmholtz equation $\Delta u + k^2 u = 0$, see for instance \cite{ruiz} for more on this connection.  As one can see, the subject of Bochner-Riesz means is a diverse one, and we will not attempt to give a comprehensive survey of it here.

\lecture{Connections with PDE estimates}

(Note: The coverage of this material here is a sketch only.  The courses of Gigliola Staffilani and Carlos Kenig will cover some of this material in far greater depth).

The restriction problem is also closely related to that of estimating solutions to linear PDE such as the wave and Schr\"odinger equations; this connection was first observed by Strichartz \cite{strichartz:restrictionquadratic}, leading to the family of estimates which bear his name.

Strichartz estimates are concerned with the spacetime $L^p$ estimates for solutions to linear PDE.  For sake of discussion we shall restrict our attention\footnote{There is an equally extensive theory for the wave equation $u_{tt} - \Delta u = 0$ and the Airy equation $u_t + u_{xxx} = 0$, and indeed one has Strichartz estimate for any dispersive equation, but we will not have the space to discuss all of these equations here.} to the \emph{Schr\"odinger equation}
$$ i u_t + \Delta u = F; \quad u(0,x) = u_0(x),$$
where the \emph{initial data} $u_0(x)$ and the \emph{forcing term} $F(t,x)$ are given complex-valued fields, which determine the complex-valued \emph{solution} $u(t,x)$.  For this discussion we shall work globally in spacetime $(t,x) \in \R \times \R^n$, although it is also of importance to consider the situation where the time variable is localized to a compact interval.

The function $u$ can be determined explicitly from $u_0$ and $F$ by several means (for instance, by Duhamel's formula); here we shall use the spacetime Fourier transform (interpreted in the distributional sense)
$$ \tilde u(\tau,\xi) := \int_\R \int_{\R^n} e^{-2\pi i (x \cdot \xi + t\tau)} u(t,x)\ dx dt.$$
Then the Schrodinger equation becomes
$$ -2\pi \tau \tilde u(\tau,\xi) - 4\pi^2 |\xi|^2 \tilde u(\tau,\xi) = \tilde F(\tau,\xi),$$
and thus in the distributional sense we have
$$ \tilde u(\tau,\xi) = - \frac{1}{2\pi} p.v. \frac{F(\tau,\xi)}{\tau - 2\pi |\xi|^2} + m(\xi) \delta(\tau - 2\pi |\xi|^2)$$
for some function $m(\xi)$.  To determine $m$, we use the initial condition
$$ \int_\R \tilde u(\tau,\xi)\ d\tau = \hat u_0(\xi).$$
and so we obtain the complete solution
$$ \tilde u(\tau,\xi) = \delta(\tau - 2\pi |\xi|^2)
(\hat u_0(\xi) + \frac{1}{2\pi}
p.v. \int_\R \frac{\tilde F(\tau',\xi)}{\tau' - 2\pi |\xi|^2}\ d\tau' - \frac{1}{2\pi} p.v. \frac{F(\tau,\xi)}{\tau - 2\pi |\xi|^2}.$$
This expression contains several recognizable components.  The contribution of third term,
$$ p.v. \frac{F(\tau,\xi)}{\tau - 2\pi |\xi|^2}$$
to $u(t,x)$ can be thought of what is more or less a Bochner-Riesz operator of order -1 (but for the paraboloid $S_{parab} := \{\tau = 2\pi |\xi|^2\}$ instead of the sphere) applied to $F$, although this is a slight oversimplification as it ignores the cancellation in the Hilbert singularity $p.v. \frac{1}{\tau}$.  The contribution of the first term,
$$ \delta(\tau - 2\pi |\xi|^2) \hat u_0(\xi)$$
to $u$ is basically the adjoint restriction operator $F \mapsto (F d\sigma)^\vee$ for the paraboloid $S_{parab}$, applied to the function $\hat u_0(\xi)$.  The contribution of the second term,
$$ \delta(\tau - 2\pi |\xi|^2)
p.v. \int_\R \frac{\tilde F(\tau',\xi)}{\tau' - 2\pi |\xi|^2}\ d\tau')$$
is almost (but not quite) the $TT^*$ of the restriction operator $u \mapsto \tilde u|_{S_{parab}}$, applied to $F$.  (This is slightly inaccurate, but becomes more truthful if the Hilbert singularity $p.v. \frac{1}{\tau}$ (the Fourier transform of the signum function) is replaced by the Dirac delta function $\delta(\tau)$ (the Fourier transform of the constant function 1).

Thus the solution operator $(u_0, F) \mapsto u$ incorporates several of the operators discussed earlier, except that when dealing with the adjoint restriction operator one must first apply a Fourier transform.  This allows one to link restriction estimates $R^*_{S_{parab}}(q' \to p')$ to estimates on the Schr\"odinger equation when $q'=2$, thanks to Plancherel's theorem; unfortunately, though, the restriction estimates when $q' \neq 2$ do not seem to have much application to these equations yet.  For simplicity we restrict our discussion to the homogeneous case when $F = 0$, so that our formula for $u$ simplifies to
$$ \tilde u(\tau,\xi) = \delta(\tau - 2\pi |\xi|^2) \hat u_0(\xi).$$
Starting with the Tomas-Stein restriction theorem $R^*_{S_{parab}}(2 \to \frac{2(n+2)}{n})$, see \eqref{ts}, we can thus obtain the \emph{Strichartz estimate}
$$ \| u \|_{L^{2(n+2)/n}(\R \times \R^n)} \lesssim \| u_0 \|_{L^2(\R^n)}.$$
This Strichartz estimate can be generalized in a number of ways. One important way is to give a separate Lebesgue exponent to the time and space exponent separately.  For instance, we have the identity
$$ \| u \|_{L^\infty_t L^2_x} = \| u_0 \|_{L^2(\R^n)}$$
just because the free Schr\"odinger equation $iu_t + \Delta u = 0$ preserves the $L^2_x$ norm.  Interpolating this estimate we can already get a certain number of mixed-norm estimates of the form
$$ \| u \|_{L^q_t L^r_x} \lesssim \| u_0 \|_{L^2(\R^n)};$$
this interpolation argument does not give the full range of estimates available, but one can modify the Tomas-Stein argument to obtain the optimal range of estimates here (see \cite{tao:keel}).  One can also combine these estimates with Sobolev embedding to obtain similar estimates when the initial data is in a Sobolev space $H^s(\R^n)$.

These estimates (and their generalizations when one restores the forcing term $F$) are particularly useful in analyzing the solutions to semi-linear Schr\"odinger equations such as
$$ i u_t + \Delta u = F(u); \quad u(0,x) = u_0(x),$$
as they often allow the effect of the non-linear term $F(u)$ to be iterated away in a suitable space such as $L^q_t L^r_x$.  This is however a huge topic in and of itself and we will not attempt to survey it here.

It may also be useful to obtain estimates where the $L^2$-type control on the initial data $u_0$ is replaced by $L^p$ control.  This may look a lot like the restriction theorems $R^*_S(p \to q)$ studied earlier, but these estimates are not the same (and in fact are substantially more difficult) because of the Fourier transform intervening between the solution operator for the Schr\"odinger or wave equation, and the corresponding adjoint restriction operator.  One conjectured estimate of this type is Sogge's \emph{local smoothing estimate} \cite{sogge:smoothing}, which addresses solutions $u(t,x)$ to the wave equation $u_{tt} - \Delta u = 0; \quad u(0,x) = f(x); \quad u_t(0,x) = 0$ in $n$ spatial dimensions, and conjectures the estimate
$$ \| u \|_{L^p([1,2] \times \R^n)} \leq C_{p,\eps} \| (1 + \sqrt{-\Delta})^\eps f \|_{L^p(\R^n)}$$
for all $\eps > 0, n(\frac{1}{2}-\frac{1}{p}) - \frac{1}{2}$ and $2 \leq p \leq \infty$.  This conjecture is easy when $p=2$ or $p =\infty$; the most interesting case is when $p = 2n/(n-1)$.  Roughly speaking, this conjecture asserts that the propagator for the wave equation is almost bounded on $L^p(\R^n)$ provided that one averages locally in time; note that without this averaging then such a claim would be highly false, as waves can certainly concentrate at a point, causing blowup of $L^p$ norms for $p>2$, while the dual situation of waves diverging from a point can cause blowup in $L^p$ norms for $p<2$.  

The local smoothing conjecture is extremely difficult; for instance, it is known to imply both the Bochner-Riesz and restriction conjectures \cite{sogge:smoothing}, \cite{tao:boch-rest}, and is also substantially stronger than the known wave equation Strichartz estimates and the circular maximal function estimates (see \cite{ss} for more on this connection); it is also related to some other operators in harmonic analysis such as convolution with curves in space \cite{oss}, \cite{tv:cone2}.  It is far from solved, even in two dimensions, although there has been some recent encouraging progress in \cite{wolff:smsub}, \cite{lw} (see also some earlier work in \cite{mock:cone}, \cite{borg:cone}, \cite{tv:cone2}).  These results are too technical to describe here in detail, but let us just mention that these results use the full panoply of techniques mentioned above (reduction to local estimates, bilinear estimates, wave packet decomposition, induction on scales, etc.).

It is also of interest to obtain estimates where the space and time norms have been reversed, e.g.
$$ \| u \|_{L^r_x L^q_t} \lesssim \| u_0 \|_{H^s(\R^n)}.$$
These estimates capture certain ``local smoothing'' properties of the solution, in that some unexpected regularity is gained in space when one averages in time.  These type of estimates are not directly related to restriction estimates, although they do have some resemblance to the maximal Bochner-Riesz estimates.  One examples of such an estimate is the \emph{Kato smoothing estimate}
$$ \| u \|_{L^\infty_x L^2_t} \lesssim \| u_0 \|_{\dot H^{-1/2}(\R)}$$
in one dimension $n=1$ (see e.g. \cite{kpv:periodic}); similar estimates hold in higher dimensions.  Another is the maximal oscillatory integral operator estimate \cite{carl}
$$ \| u \|_{L^4_x L^\infty_t} \lesssim \| u_0 \|_{\dot H^{1/4}(\R)};$$
the sharp higher-dimensional analogue of this estimate is not yet fully understood (see \cite{tv:cone2} for some recent results on this problem); interestingly, the bilinear restriction estimates for the paraboloid mentioned earlier are somewhat effective in addressing the problem, although they do not seem strong enough to provide the optimal estimates.

Just as linear restriction theorems have bilinear counterparts, so too do Strichartz estimates, where one now has two solutions $u$, $v$ to (say) the free Schr\"odinger equation 
$$i u_t + \Delta u = i v_t + \Delta v = 0$$
and one now estimates bilinear expressions such as $uv$, or perhaps a more fancy version of this expression involving differential operators (for instance, in the context of the wave equation one often wishes to estimate \emph{null forms}, which roughly speaking are those bilinear forms which vanish when tested against pairs of parallel plane waves).  One might also wish to impose some transversality condition on the supports of the Fourier transforms of the initial data $u(0)$ and $v(0)$.  These types of estimates have been extensively studied, especially when one measures the bilinear expression in a space based on $L^2_{t,x}$, for instance an $L^2$ Sobolev space\footnote{In the study of the non-linear Schr\"odinger equation it turns out that the natural spaces to use are not just the standard Sobolev spaces $H^s$ - which are based on the powers of the Laplacian - but also a variant space, called $X^{s,b}$ or $H^{s,\theta}$ in the literature, which are based on powers of the Laplacian and of the Schr\"odinger operator $i\partial_t + \Delta$.  See for instance \cite{borg:xsb}, \cite{klainerman:nulllocal}, \cite{kpv:kdv}, \cite{ginibre:survey}, \cite{tao:xsb}.}.  A typical such estimate is
$$ \| uv \|_{L^2_t H^s_x} \lesssim \| u(0) \|_{L^2_x} \| v(0) \|_{H^s_x}$$
for all $s > 1/2$ and all solutions $u,v$ to the free Schr\"odinger equation in two dimensions (see \cite{borg:nls}).  These sorts of estimates have many applications to nonlinear Schr\"odinger equations, for instance they are useful in demonstrating that the low frequencies and high frequencies of a nonlinear Schr\"odinger wave have very little interaction, and thus evolve more or less independently.

The $L^p$ bilinear restriction estimates for the paraboloid discussed in previous lectures can lead to bilinear $L^p$ estimates for solutions to the Schr\"odinger equation.  For instance, the bilinear estimate in Theorem \ref{tao:thm}
leads directly to the bilinear estimate
$$ \| u v \|_{L^p_{t,x}} \leq C_{S_1,S_2,p} \| u(0) \|_{L^2_x} \| v(0) \|_{L^2_x}$$
for all solutions $u,v$ to the free Schr\"odinger equation whose Fourier transforms are supported on disjoint compact sets $S_1, S_2$, and for $p > \frac{n+3}{n+1}$.  There is a similar result for products of solutions of the wave equation arising from Theorem \ref{wolff:thm}, see for instance \cite{tao:cone} for a discussion.  These types of estimates may well have significant 
applications to nonlinear Schr\"odinger equations in the future; at present these $L^p$ bilinear estimates have only found a few specialized applications (see e.g. \cite{planchon}).

Recently, there has been a lot of interest in developing Strichartz estimates (and their bilinear counterparts) on manifolds or when there is a potential term in the equation, both for the Schr\"odinger and wave equations; this is a vast and highly active field which is impossible to summarize here, but a small (and very incomplete) sample of recent results is \cite{rs}, \cite{st}, \cite{tataru}, \cite{kr}, \cite{smt}, \cite{bgt}, \cite{bpss}; see also the lectures of Carlos Kenig.  In such situations the spacetime Fourier transform is no longer so useful (although the temporal Fourier transform still has some uses, and the spatial Fourier transform can to some extent be replaced by the spectral decomposition of the Laplacian or Hamiltonian), and so restriction theory \emph{per se} no longer applies; however many of the techniques used in that theory still do apply, notably the $TT^*$ method and wave packet decompositions.  Good control on the fundamental solution (or of an approximate parametrix for the fundamental solution) becomes important (this is the PDE analogue of the quantity $\widehat{d\sigma}$ which played a key role in the restriction theory).  I believe that the continuing research in this area will not only generalize the Strichartz estimates, but also deepen our insight into exactly why they arise and what their range of applicability is.

Another facet of restriction theory on manifolds arises in the study of eigenfunctions of the Laplacian.  To motivate this, we return to the local restriction estimates \eqref{eq:blur-adjoint} for the sphere $S_{sphere}$ and dilate it by $R$.  This converts the local restriction estimate $R^*_{S_{sphere}}(q' \to p')$ to the estimate
$$
\| G^\vee \|_{L^{p'}(B(x_0,1))} \lesssim R^{\alpha + \frac{n-1}{q} - \frac{n}{p'}} \| G \|_{L^{q'}(N_{1}(R S_{sphere}))}.$$
In particular, if $q = 2$ and $G$ is the Fourier transform of some function $u$, then we see that we have the estimate
\begin{equation}\label{eq:u-est}
 \| u \|_{L^{p'}(B(0,1))} \lesssim R^{\alpha + \frac{n-1}{2} - \frac{n}{p'}}
\| u \|_{L^2(\R^n)}
\end{equation}
for any function $u$ whose Fourier transform is supported on the annulus $\{ R-1 \leq |\xi| \leq R+1 \}$.  This last condition is equivalent to $u$ lying in the range of the spectral projection $\chi_{[R-1,R+1]}(\sqrt{-\Delta})$ (we ignore factors of $2\pi$ for this heuristic discussion).

We can view \eqref{eq:u-est}, when $p'$ is larger than $2$, as a way to control the large values of an $L^2$-normalized function $u$ in the range of the above-mentioned projection.  In particular, this should give us some local control on the large values of eigenfunctions $\sqrt{-\Delta} u = R u$, or equivalently $\Delta u + R^2 u = 0$.

It is of some importance to extend these sorts of estimates to more general Riemannian manifolds $(M,g)$; for sake of argument suppose that $M$ is smooth and compact without boundary.  The Laplacian $\Delta$ will then of course be replaced by the Laplace-Beltrami operator $\Delta_g$, and the objective now is to obtain estimates of the form
$$ \| u \|_{L^{p'}(M)} \lesssim R^\beta \| u \|_{L^2(M)}$$
for various exponents $p$, $\beta$, where $u$ is any function in the range of the spectral projection $\chi_{[R-1,R+1]}(\sqrt{-\Delta_g})$.  The significance of such estimates is that they give some control on the distribution of eigenfunctions of the Laplace-Beltrami operator; this has relevance to some difficult problems concerning quantum ergodicity, and also to questions in analytic number theory (for instance, in estimating Hecke forms, which are speical eigenfunctions of the Laplacian on cusp domains).  These estimates are also relevant to the problem of obtaining Strichartz estimates on manifolds, which was mentioned earlier.

It is possible to make some headway on these problems using restriction techniques.  One key point is to decompose (a smoothed out version of) $\chi_{[R-1,R+1]}(\sqrt{-\Delta_g})$ - which plays a role similar to the dyadic component $m_j(D)$ of the Bochner-Riesz operator in the Euclidean theory - as a combination of wave propagators $\exp(\pm i t \sqrt{-\Delta g})$ for short times $t = O(1)$, and then analyze the fundamental solution to those wave operators (or to a parametrix thereof); see e.g. \cite{sogge:Fourier}.  This gives good results for special manifolds such as the sphere, in which the spectrum of the Laplacian is highly concentrated on a sparse set.  For more general manifolds, the theory is less satisfactory; it seems that we have to localize the spectral projections to a narrower interval than $[R-1,R+1]$, but this then seemingly requires us to analyze the wave equation for much longer times $t \gg 1$, which is largely beyond the ability of our current technology.  A breakthrough in this direction would have significant impact on the above-mentioned fields, although it seems that this requires techniques quite distinct from those used in the existing restriction theory.  We should warn however that we do not expect all the results from the Euclidean theory to carry over to curved space, see for instance \cite{sogge:nikodym} for some counterexamples.

\end{document}